\documentclass[a4paper,reqno,oneside]{amsart}
\usepackage{amsmath,amsrefs,geometry}
\usepackage{graphicx}
\usepackage{mathrsfs}
\usepackage{amssymb,amsmath, amsfonts, amsthm}
\usepackage{latexsym}
\usepackage{color}
\usepackage[dvips]{epsfig}
\usepackage{graphicx}

\allowdisplaybreaks[1]

\numberwithin{equation}{section}

\renewcommand{\(}{\left(}
\renewcommand{\)}{\right)}
\renewcommand{\[}{\left[}
\renewcommand{\]}{\right]}

\newtheorem{theorem}{Theorem}[section]
\newtheorem{proposition}[theorem]{Proposition}
\newtheorem{lemma}[theorem]{Lemma}
\newtheorem{remark}[theorem]{Remark}

\renewcommand{\le}{\leqslant}
\renewcommand{\ge}{\geqslant}

\renewcommand{\S}{\Sigma}

\newcommand{\A}{{\mathcal A}}

\newcommand{\tm}{\tilde{\mu}}

\newcommand{\te}{\tilde{e}}
\newcommand{\D}{\mathcal{D}}
\newcommand{\M}{\mathcal{M}}

\renewcommand{\S}{\mathcal{S}}

\newcommand{\beq}{\begin{equation}}
\newcommand{\eeq}{\end{equation}}
\newcommand{\beqs}{\begin{equation*}}
\newcommand{\eeqs}{\end{equation*}}
\newcommand{\beqn}{\begin{eqnarray}}
\newcommand{\eeqn}{\end{eqnarray}}
\newcommand{\beqns}{\begin{eqnarray*}}
\newcommand{\eeqns}{\end{eqnarray*}}
\newcommand{\bdoc}{\begin{document}}
\newcommand{\edoc}{\end{document}}
\newcommand{\be}{\begin{enumerate}}
\newcommand{\ee}{\end{enumerate}}
\newcommand{\bdescr}{\begin{description}}
\newcommand{\edescr}{\end{description}}
\newcommand{\ba}{\begin{array}}
\newcommand{\ea}{\end{array}}
\newcommand{\intR}{\int_{\mathbb R^N}}
\newcommand{\R}{\mathbb R^N}
\newcommand{\e}{\epsilon}
 \renewcommand{\(}{\left(}
\renewcommand{\)}{\right)}
\renewcommand{\[}{\left[}
\renewcommand{\]}{\right]}

\begin{document}
\title[ Bubbling solutions for supercritical problems on manifolds]{Bubbling solutions  for supercritical problems on manifolds\footnotemark}
\author{Juan D\'avila}
\address[Juan D\'avila]{Departamento de Ingenier\'ia Matem\'atica and Centro de Modelamiento Matem\'atico, Universidad de Chile, Casilla 170 Correo 3, Santiago, Chile}
\email{jdavila@dim.uchile.cl}
\author{Angela Pistoia}
\address[Angela Pistoia] {Dipartimento SBAI, Universt\`{a} di Roma ``La Sapienza", via Antonio Scarpa 16, 00161 Roma, Italy}
\email{pistoia@dmmm.uniroma1.it}
\author{Giusi Vaira}
\address[Giusi Vaira] {Dipartimento di Matematica ``G. Castelnuovo", Universit\`{a} di Roma ``La Sapienza", Piazzale A. Moro 1, 00161 Roma, Italy}
\email{vaira@mat.uniroma1.it}

\begin{abstract}
Let $(\M,g)$  be a $n-$dimensional compact Riemannian manifold without boundary and $\Gamma$ be a non degenerate closed geodesic of $(\M,g)$.
We prove that    the  supercritical problem
$$-\Delta _gu+h u=u^{\frac{n+1}{n-3}\pm\e},\ u>0,\ \hbox{in}\ (\M,g)$$
has a solution that concentrates along $\Gamma$ as $\e$ goes to zero, provided the function $h$  and the sectional curvatures
along $\Gamma$ satisfy a suitable condition.
A connection  with the solution of a class of periodic O.D.E.'s with singularity of attractive or repulsive type is established.

 \end{abstract}


\date{\today}

 \maketitle

\footnotetext{The first author was supported by Fondecyt grant 1130360 and Fondo Basal CMM. The second and the third authors have been partially supported by the Gruppo Nazionale per l'Analisi Matematica, la Probabilit\'a  e le loro Applicazioni (GNAMPA) of the Istituto Nazionale di Alta Matematica (INdAM).
}

\maketitle

\section{Introduction and statement of main results}\label{intro}

We deal with the semilinear elliptic equation
\begin{equation}\label{p}
 -\Delta _gu+h u=u^{p-1},\ u>0,\ \hbox{in}\ (\M,g)
 \end{equation}
where   $(\M,g)$  is a $n-$dimensional compact Riemannian manifold without boundary, $h$ is a $C^1-$real function on $\M$   such that $-\Delta_g+h$  is coercive and $p>2$.

For any $p\in(2,2^*_n)$, where $2^*_n:={2n\over n-2}$ if  $n\ge3$ and $ 2^*_n:=+\infty$ if $n=2$,   problem \eqref{p}  has a solution, which can be found by minimization of
$$
\mathcal I_p(u)={
\int\limits_\M\(|\nabla_g u|^2+h u^2\)d\sigma_g\over \(\int\limits_\M | u|^pd\sigma_g\)^{2/p}}
$$
over $H^1_g(\M)\setminus\{0\}$, using the compactness of the embedding $H^1_g(\M)\hookrightarrow L^p_g(\M)$.

In the critical case, i.e.  $p=2^*_n,  $  the situation turns out to be more delicate. In particular, the existence of solutions is related to the position of the potential $h $ with respect to the geometric potential $h_g:= {m-2\over 4(m-1)} R_g$, where $R_g$ is the scalar curvature of the manifold.

{
If $h\equiv h_g$, then problem \eqref{p}  is referred to as  the   Yamabe problem \cite{Yam} and it has always a solution. After Trudinger \cite{Tru} discovered a gap in the argument in \cite{Yam} and gave a proof under some conditions on $(\M,g)$, Aubin \cite{Aub1,Aub2} showed that whenever $Q(\M,g)<Q(S^n,g_0)$, where $(S^n,g_0)$ is the standard sphere and
$$
Q(\M,g) := \inf_{u\in H^1_g(\M)\setminus\{0\}} I_{2^*_n}(u) ,
$$
there is a solution to the problem, and proved that this holds if $n\geq 6 $ and $(\M,g)$ is not locally conformally flat. Finally, Schoen \cite{Sch1} gave a proof in full generality using the Positive Mass Theorem \cite{schoen-yau}.
}

When $h< h_g$ somewhere in $M$, existence of a solution is guaranteed by a minimization argument, arguing as in  Aubin \cite{Aub1,Aub2}. The situation is extremely delicate when  $h \ge h_g$ everywhere in $\M$, because blow-up phenomena can occur as pointed out by Druet in \cite{D1,D2}.

The supercritical case $p>2^*_n$ is even more difficult to deal with.
A first result in this direction is a perturbative result   due to Micheletti, Pistoia and V\'etois \cite{mpv}.
They consider the almost critical problem \eqref{p} when
$p=2^*_n \pm\epsilon$, i.e.
if $p=2^*_n-\epsilon$ the problem \eqref{p} is slightly subcritical and if $p=2^*_n+\epsilon$ the problem \eqref{p} is slightly supercritical. They prove the following results.
\begin{theorem}\label{mipive}[Micheletti, Pistoia and V\'etois \cite{mpv}]
 Assume
 $n\ge6  $ and  $\xi_0\in M$ is a non degenerate critical point of $h-{m-2\over 4m}R_g.$
Then
 \begin{itemize}
 \item[(i)] if $ h(\xi_0)>{n-2\over 4n}R_g(\xi_0)$ then   the slightly subcritical   problem \eqref{p} with $p=2^*_{n}-1-\e,$ has a solutions $u_\e$ which concentrates at $\xi_0$ as $\e\rightarrow 0$,
 \item[(ii)] if $ h(\xi_0)<{n-2\over 4n}R_g(\xi_0)$ then the slightly supercritical problem \eqref{p} with $p=2^*_{n}-1-\e,$ has a solutions $u_\e$ which concentrates at $\xi_0$ as $\e\rightarrow 0$.
\end{itemize}
\end{theorem}

Now, for any integer $0\le k\le n-3 $ let $2^*_{n,k}={2(n-k)\over n-k-2} $  be the \textit{$(k+1)-$st critical exponent}.
We remark that $2^*_{n,k}=2^*_{n-k,0}$   is nothing but the critical exponent for the Sobolev embedding $H^1_h(\mathcal N) \hookrightarrow L^q_h(\mathcal N)$ in a compact $(n-k)-$dimensional Riemannian manifold $ (\mathcal N, h)$.
In particular,  $ 2^*_{n,0}={2 n\over n- 2}$    is the usual Sobolev critical exponent.\\
We can summarize the results proved by Micheletti, Pistoia and V\'etois just saying that  problem \eqref{p} when $p\to 2^*_{n,0}$ (i.e. $k=0$) has positive  solutions   blowing-up at points.
    Note that   a  point  is a   $0-$dimensional manifold!\\\\  A natural question arises:\\
  \textit{ does  problem \eqref{p} have solutions blowing-up at  $k-$dimensional submanifolds   when $p\to 2^*_{n,k}$?}\\\\

In the present paper,   we give a positive answer when $k=1.$ More precisely, we prove that if $p\to 2^*_{n,1}$ problem
\eqref{p} has a solution which concentrates along a geodesic $\Gamma$ of the manifold provided  $h$  satisfies a suitable  condition.
 Let us state our main result.\\

{
We consider the problem \eqref{p} with $p=2^*_{n,1}\pm\e$, i.e. }
 \begin{equation}\label{pb}
 -\Delta _gu+h u=u^{\frac{n+1}{n-3}\pm\e},\ u>0 \ \hbox{in}\ (\M,g)
 \end{equation}
We will say that problem \eqref{pb} is   slightly \textit{$2nd-$supercritical}   if $p=2^*_{n,1}+\e$
and it is  slightly \textit{$2nd-$subcritical}   if $p=2^*_{n,1}-\e$.\\

Let $\Gamma $ be a closed nontrivial { simple} geodesic in $\M.$    We define the function
\begin{equation}\label{gamma0}
 \sigma (x_0)=
   h(x_0) -{(n-3)(n-4)\over 3(n-2)}\(  R_{g}(x_0) -{N\over 4}  Ric(\dot \gamma(x_0),\dot \gamma(x_0))\)
\end{equation}
where     $R_g$ is the scalar curvature and {  $Ric$ denotes the Ricci tensor.}


%

  Let    $a_n:={8(n-2)\over (n-3)(n+1)}
$ and $b_n:={ (n-3)^2(n-5)\over 4(n+1)} $ (see \eqref{costanti} and Remark \eqref{gamma}). \\
  We   introduce the   periodic ODE problem
 \begin{equation}\label{ip-meno}
 \left\{\begin{aligned}
&-\ddot \mu+a_n\sigma  \mu- {b_n\over\mu}=0\quad \hbox{in}\ [0,2\ell]  \\
&\mu>0 \quad \hbox{in}\ [0,2\ell]\\
&\mu(0)=\mu(2\ell),\ \dot\mu(0)=\dot\mu(2\ell)\\
\end{aligned}\right.
\end{equation}
which has a \textit{ singularity of attractive type} at the origin and the  periodic ODE problem
\begin{equation}\label{ip-piu}
 \left\{\begin{aligned}
&-\ddot \mu+a_n\sigma  \mu+ {b_n\over\mu}=0\quad \hbox{in}\ [0,2\ell] \\
&\mu>0 \quad \hbox{in}\ [0,2\ell]\\
&\mu(0)=\mu(2\ell),\ \dot\mu(0)=\dot\mu(2\ell)\\
\end{aligned}\right.
\end{equation}
which has a \textit{singularity of repulsive type} at the origin.\\

Solvability of the slightly \textit{$2nd-$subcritical} problem is strictly related with solvability  of  \eqref{ip-meno} with \textit{attractive   singularity}, while solvability of the slightly \textit{$2nd-$supercritical} problem
is strictly related with solvability of \eqref{ip-piu} with \textit{repulsive   singularity}.\\

As usual in this kind of problem, we also need to  assume    a gap condition of the form
\beq\label{nonresonanza}
|\e k^2-\kappa^2|>\nu \sqrt{\e},\qquad k=1, 2, \ldots \ldots
\eeq
where $\kappa>0$ is given explicitly in Lemma \ref{l0} and $\nu$ is positive.\\

Now we can state our main result.
\begin{theorem}\label{principale}
Let $n\geq 8.$  Let
 $\Gamma $ be a simple closed, non degenerate geodesic of $\M$ (see \eqref{jacobi}).
\begin{itemize}
\item[(i)] Assume the problem   \eqref{ip-meno}
has a non degenerate positive solution $\mu_0$.
 Then, for any $\nu>0$ there exists $\epsilon_0>0$ such that for any $\e\in(0,\e_0)$ which    satisfies condition \eqref{nonresonanza},  the slightly \textit{$2nd-$subcritical} problem \eqref{pb} with $p=2^*_{n,1}-1-\e,$ has a solution $u_\e$ that concentrates along $\Gamma$ as $\e\rightarrow 0$.
\item[(ii)] Assume the problem   \eqref{ip-piu}
has a non degenerate positive solution $\mu_0$.
Then, for any $\nu>0$ there exists $\epsilon_0>0$ such that for any $\e\in(0,\e_0)$ which    satisfies condition \eqref{nonresonanza}, the  slightly \textit{$2nd-$supercritical} problem \eqref{pb} with $p=2^*_{n,1}-1+\e,$ has a solution $u_\e$ that concentrates along $\Gamma$ as $\e\rightarrow 0$.
\end{itemize}

 \end{theorem}
Moreover, the solution $u_\e$  can be described in Fermi  coordinates  as follows.
 Given $\xi\in \Gamma$ there is a natural splitting $T_\xi M=T_\xi \Gamma\oplus N_\xi\Gamma$ into the tangent and normal bundle over $\Gamma.$
It is useful
to introduce a local system of coordinates    near $\Gamma.$
Let $\gamma:[0, 2\ell]\rightarrow \M$ be an arclenght parametrization of $\Gamma$, where $2\ell$ is the lenght of $\Gamma$. We denote by $E_0$ a unit tangent vector to $\Gamma$. In a neighborhood of a point $\xi$ of $\Gamma$ we give an orthonormal basis $E_1,\dots,E_N$ of $N_q\Gamma.$ We can assume that the $E_i$'s are parallel along $\Gamma,$ i.e. $  \nabla_{E_0}E_i=0$ for any $i=1,\dots,N.$ The geodesic condition for $\Gamma$ translates into the condition $  \nabla_{E_0}E_0=0.$ Here $\nabla$ is the connection associated with the metric $g.$
Moreover,  the non degeneracy of $\Gamma$ is equivalent to say that  the linear equation
\begin{equation}
\label{jacobi}
\mathcal J\phi:=\nabla^2_{E_0}\phi+R(\phi,E_0)E_0=0\ \hbox{has only the trivial solution on all of $\Gamma.$}
\end{equation}
 Here $\mathcal J$ is the Jacobi operator on $\Gamma$ corresponding to the second variation of the length  functional on curves. For a generic metric $g$ on $\M$ it is well known that all closed geodesics are non degenerate.\\
To parametrize a neighborhood of a point of $\Gamma$ in $M$ we
  define the  \textit{ Fermi coordinates}
\begin{equation}
\label{fermi0}
F(x_0,x_1,\ldots, x_N)={\rm exp}_{\gamma(x_0)}\left(\sum_{i=1}^N x_i E_i(x_0)\right)
\end{equation}
where ${\rm exp}_{\gamma(x_0)}$ is the exponential map in $\M$ through the point $\gamma(x_0)$.\\\\

The solution constructed in Theorem~\ref{principale} has the expansion
\begin{equation*}
u_\e(x_0, x)=\mu_\e^{-\frac{N-2}{2}}w\left(\mu_\e^{-1}(x-d_\e)\right)+o(1),\end{equation*}
where $$\mu_\e(x_0)\sim \sqrt{\e}\mu_0(x_0)\ \hbox{and}\ d_{\e _k}(x_0)\sim \e  {d}_k(x_0), k=1, \ldots, N \qquad $$ where $\mu_0$ solves  either problem
\eqref{ip-meno} in the slightly \textit{$2nd-$subcritical} case or problem \eqref{ip-piu} in the slightly \textit{$2nd-$supercritical} case, the
$ {d}_j$'s   are smooth functions of $x_0$ and $w$ is the   standard bubble
\begin{equation}\label{stanbu}
w(y)=c_N\frac{ 1}{\left(1+|y|^2\right)^{\frac{N-2}{2}}},\qquad y\in \mathbb R^N,\qquad c_N=[N(N-2)]^{\frac{N-2}{4}},\end{equation} which is the radial solution of the critical problem $\Delta w+w^p=0$ in $\mathbb R^N.$ We point out  that $N=n-1.$\\

Since the existence of   solutions to   singular problems  \eqref{ip-meno} or \eqref{ip-piu} plays a crucial role in the construction of the solution, in particular in the choice of the concentration parameter
$\mu_\e$, it is important to  point out  that  existence of solutions to problems  \eqref{ip-meno} or \eqref{ip-piu} is strictly linked with the sign of the function $\sigma $ defined in \eqref{gamma0}, as it is showed  in the following Theorem, whose proof is given in Section \ref{ode}.

\begin{theorem}\label{main2}

If
\begin{equation*}
 \min\limits_{x_0\in\mathbb R}\sigma (x_0) >0,
\end{equation*}
then problem   \eqref{ip-meno} has a non degenerate solution.\\

If $h^*\in C^2(M)$ is such that
\begin{equation*}
-\({(k+1)\pi\over2\ell}\)^2<\min\limits_{x_0\in\mathbb R}\sigma_{h^*}(x_0) \le \max\limits_{t\in\mathbb R}\sigma_{h^*}(x_0)< -\({k\pi\over2\ell}\)^2 <0, \end{equation*}
 then  for most functions $h\in C^2(M)$ with  $  \|h-h^*\|_{C^0(M)}\le r,$ provided $r$ is small enough,   the  problem   \eqref{ip-piu} has a   non degenerate solution.

\end{theorem}

   As far as we know, Theorem \ref{principale} is the first result about existence of solutions to \eqref{p} which concentrate along geodesic of  the manifold $M$ when the exponent $p$ approaches the $2nd-$critical exponent from above. Indeed,    in the Euclidean setting,  del Pino, Musso and Pacard in \cite{DMP}  built bubbling solutions for a Dirichlet problem when the exponent is close to but less than   the second critical exponent.
{ Solutions concentrating in higher dimensional sets and the gap condition have been found in elliptic problems in the Euclidean setting. We mention among, among many results, \cite{malchiodi-montenegro-cpam,malchiodi-montenegro-duke,mahmoudi-malchiodi-adv,malchiodi-gafa} for a Neumann singular perturbation problem and \cite{delpino-kowalczyk-wei-cpam}
for a Sch\"odinger equation in the plane.}
\\\\
We conjecture that our result can be extended to higher  $k-$dimensional minimal submanifold of $\M.$ Indeed,  arguments developed by Del Pino, Mamhoudi and Musso in \cite{DMM} in the Euclidean setting for a Neumann problem could also be applied to equation \eqref{p}.
\\\\
The proof of our result  relies on the    infinite-dimensional reduction   developed by Del Pino, Musso and Pacard in \cite{DMP}.
We omit many details of proof, because they   can be found, up to some minor modifications, in \cite{DMP}. We only compute what can not be
deduced by known results.\\

The paper is organized as follows.
In Section \ref{ode} we study the singular problems
\eqref{ip-meno} and \eqref{ip-piu}. In Section \ref{approssimante} we build the approximate solution close to the geodesic
 and in Section \ref{sec-errore} we estimate the error.
 Then, in Section \ref{gluing} we reduce the problem to a suitable infinite dimensional set of parameters and in Section \ref{equazioni} we study the reduced problem. Section \eqref{lineare} is devoted to the study of a linear problem.\\\\

{\bf Notation}
\begin{itemize}
\item
For sums we use the standard convention of summing terms where repeated indices appear.
\item We will denote by   $L^\infty_{2\ell}(\mathbb R),$ $C^0_{2\ell}(\mathbb R)$ and $C^2_{2\ell}(\mathbb R)$   the Banach space of $2\ell-$periodic $L^\infty$, $C^0$ and $C^2$ functions, respectively. We will set $\|u\|_\infty:=\sup\limits_{\mathbb R}|u|,$  for any   $2\ell-$periodic bounded  function $u$.
 \end{itemize}

\section{A periodic ODE   with repulsive or attractive singularity}\label{ode}
Let us consider the periodic boundary value problem
\begin{equation}\label{ode1}
\left\{\begin{aligned}
&-\ddot \mu+\sigma \mu-{c\over\mu}=0 \quad \hbox{in}\ [0,2\ell]\\
& \mu>0 \quad \hbox{in}\ [0,2\ell]\\
&\mu(0)=\mu(2\ell),\ \dot\mu(0)=\dot\mu(2\ell)\\
\end{aligned}\right.
\end{equation}
where $c\in \mathbb R$ and $\sigma\in C^0_{2\ell}(\mathbb R).$
The following existence result holds true.
\begin{proposition}\label{habets+delpino}
Assume either
\begin{equation}\label{piu}
 \min\limits_{t\in\mathbb R}\sigma(t) >0\ \hbox{and}\  c>0
\end{equation}
or
\begin{equation}\label{meno}
-\({(k+1)\pi\over2\ell}\)^2<\min\limits_{t\in\mathbb R}\sigma(t) \le \max\limits_{t\in\mathbb R}\sigma (t)< -\({k\pi\over2\ell}\)^2 <0\ \hbox{and}\  c<0
\end{equation}
for some integer $k$.
Then problem \eqref{ode1} has a periodic solution $\mu_0\in C^2_{2\ell}(\mathbb R).$  \\
\end{proposition}
\begin{proof}
If \eqref{piu} holds,   the claim follows   by Proposition 1 of \cite{HS} and if \eqref{meno} holds  the claim  follows   by Theorem 1.1 of \cite{DMaMo}.
\end{proof}
Let us consider the linearization of problem \eqref{ode1} around $\mu_0,$   namely the linear periodic boundary value problem
\begin{equation}\label{ode2}
\left\{\begin{aligned}
&-\ddot \mu+\(\sigma  +{c\over\mu_0^2}\)\mu =0\quad \hbox{in}\ [0,2\ell]\\
&\mu(0)=\mu(2\ell),\ \dot\mu(0)=\dot\mu(2\ell)\\
\end{aligned}\right.
\end{equation}
The solution $\mu_0$ is non degenerate if and only if the problem \eqref{ode2} has only the trivial solution.

\begin{proposition}\label{nondegene}
\begin{itemize}
\item [(i)]
If   \eqref{piu} holds,  then the solution $\mu_0$ is non degenerate.
\item[(ii)] Let $\sigma ^*\in C^0_{2\ell}(\mathbb R)$ and $c\in\mathbb R$ as in \eqref{meno}.
The set
$$ \left\{\sigma \in  B(\sigma ^*,r)\ :\  \hbox{ all the positive solutions of \eqref{ode1} are nondegenerate}
 \right\}$$
is a dense subset of the ball $B(\sigma ^*,r):=\left\{\sigma \in C^0_{2\ell}(\mathbb R)\ :\     \|\sigma  -\sigma ^*\|_\infty\le r \right\}$  provided the radius $r$ is small enough.
\end{itemize}
 \end{proposition}
\begin{proof}
 (i) follows immediately by the maximum principle.\\
 Let us prove (ii). We shall use the following abstract transversality theorem previously used by Quinn \cite{Q}, Saut and Temam \cite{ST} and Uhlenbeck \cite{U}.
  \begin{theorem}\label{tran}
 Let $X,Y,Z$ be three Banach spaces and $U\subset X,$ $V\subset Y$ open subsets.
 Let $F:U\times V\to Z$ be a $C^\alpha-$map with $\alpha\ge1.$ Assume that

 \begin{itemize}
 \item[$(\iota)$] for any $y\in V$, $F(\cdot,y):U\to Z$ is a Fredholm map of index $l$ with $l\le\alpha;$
 \item[$(\iota\iota)$] $0$ is a regular value of $F$, i.e. the operator $F'(x_0,y_0):X\times Y\to Z$ is onto at any point $(x_0,y_0)$ such that $F(x_0,y_0)=0;$
 \item[$(\iota\iota\iota)$] the map  $\pi\circ i:F^{-1}(0)\to Y$ is $\sigma-$proper, i.e.  $F^{-1}(0)=\cup_{\eta=1}^{+\infty} C_\eta$
 where $C_\eta$ is a closed set and the restriction $\pi\circ i_{|_{C_\eta}}$ is proper for any $\eta$; here $i:F^{-1}(0)\to Y$ is the canonical embedding and $\pi:X\times Y\to Y$ is the projection.
 \end{itemize}

Then the set
$\Theta:=\left\{y\in V\ :\ 0\ \hbox{is  a regular value of } F(\cdot,y)\right\}$
 is a  residual subset of $V$, i.e. $V\setminus \Theta$ is a countable union of closet subsets without interior points.

\end{theorem}

In our case the $C^2-$ function  $F$ is defined by
$$F:C^2_{2\ell}(\mathbb R)\times C^0_{2\ell}(\mathbb R) \to C^0_{2\ell}(\mathbb R),\quad F( \mu,\sigma  ):= -\ddot \mu+ \sigma \mu -{c\over\mu },
$$
$X =C^2_{2\ell}(\mathbb R)$ and $U=\{\mu\in C^2_{2\ell}(\mathbb R)\ :\ \min_{\mathbb R}\mu>0\},$  $Y=Z= C^0_{2\ell}(\mathbb R)$  and $V=B(\sigma ^*,r)$
where   $r$ is small enough so that condition \eqref{meno} holds for any $\sigma \in V.$   \\
It is not difficult to check that   for any $\sigma \in V$ the map $\mu\to F(\mu,\sigma )$ is a Fredholm map of index $0$ and then assumption $(\iota)$ holds.
 Let us prove assumption $(\iota\iota)$. We fix $(\mu_0,\sigma _0)\in U\times V$ such that $F(\mu_0,\sigma _0)=0.$ The derivative $D_\sigma  F(\mu_0,\sigma _0): C^0_{2\ell}(\mathbb R) \to C^0_{2\ell}(\mathbb R)$ is the linear map defined by
 $D_\sigma  F(\mu_0,\sigma _0)[\sigma ]=\sigma  \mu_0 $
and it is   surjective, because $ \mu_0>0 .$\\
As far as it concerns assumption  $(\iota\iota\iota)$, we have that
$$F^{-1}(0)= \cup_{m=1}^{+\infty} \left\{\(C_m\times B_m\)\cap F^{-1}(0)\right\} $$
where
$$C_m =  \left\{\mu\in C^2_{2\ell}(\mathbb R)\ :\ {1\over m}\le \min_{\mathbb R}\mu \le\max_{\mathbb R}\mu  \le m\right\}\ \hbox{and}\ B_m=  \overline  { B  \(\sigma ^*,r-{1\over m} \)}  .$$
We can show  that the restriction
$\pi\circ i_{|_{C_m}}$ is proper, namely if the sequence $(\sigma _n)\subset  B_m $ converges to $\sigma $ and the sequence $(\mu_n)\subset  C_m$ is such that
$F(\mu_n,\sigma _n)=0$ then there exists a subsequence of $(\mu_n)$ which converges to $\mu\in C_m $ and $F(\mu ,\sigma )=0.$\\
That concludes the proof.
\end{proof}

\begin{proof}[Proof of Theorem \ref{main2}]
It follows immediately by Proposition \ref{habets+delpino} and Proposition \ref{nondegene}.
\end{proof}

\section{Construction of the approximate solution close to the geodesic}\label{approssimante}

This section is devoted to the construction of an approximation for a solution to the problem \eqref{pb} in a neighborhood of the geodesic.\\

 \subsection{The problem near to the geodesic}

Let us consider the system of Fermi coordinates $(x_0, x )$ introduced in \eqref{fermi0}. In this language the geodesic $\Gamma$ is represented by the $x_0-$ axis. We recall that $x_0$ denotes the arclenght of the curve, $2\ell$ represent the total length of the geodesic and $x=(x_1,\dots,x_N)\in \mathbb R^N.$
Let us introduce a neighborhood of the geodesic $\Gamma$ in this system of coordinates
\begin{equation}\label{deta}
D :=\left\{(x_0, x)\in \mathbb R\times\mathbb R^N\ :\ x_0\in [-\ell, \ell],\  |x| <\hat\delta \right\},\end{equation}
where $\hat\delta>0$ is a fixed small number.
Then for a function defined in $D $  we write
$$\tilde{u}(x_0, x)=u(F(x_0, x)) $$ and  we extend $\tilde{u}$ in a { satisfying the following periodicity condition
$$
\tilde u(2\ell,x) = \tilde u(0,A x)
$$
where $A=(a_{ij})$ is the invertible matrix defined by the requirement
\begin{equation}
\label{defA}
E_i(2\ell) = \sum_{j=1}^N a_{ji} E_j(0).
\end{equation}
}

Therefore, if $u$ solves
equation \eqref{pb}   in the neighborhood $D $ of the geodesic, then $\tilde u$ solves
\beq\label{3pal}
\left\{
\begin{array}{lr}
\partial_{00} \tilde u +\Delta_x\tilde u+ B(\tilde u)-h \tilde u +f_\e(\tilde{u}) =0\  \mbox{in}\ D \\
\tilde u(x_0+2\ell, x)= \tilde u(x_0, A x)
\  \mbox{for any}\  (x_0, x)\in D
\end{array}
\right.
\eeq
where   $f_\e(s):=(s^+)^{p\pm\e}.$
For the sake of simplicity, we will refer to  $f_\e(s):=(s^+)^{p+\e}$ as the \textit{supercritical case} and to $f_\e(s):=(s^+)^{p-\e}$ as the \textsl{subcritical case}.\\
In \eqref{3pal} $B$ is a second order linear operator defined in the following Lemma
 \begin{lemma}\label{fermi}
Let $u$ be a smooth function. Then for any $(x_0, x)\in D$ we have
\begin{align*}
  \Delta_g u &= \partial_{00} \tilde u +\Delta_x\tilde u +B(\tilde u)
  \end{align*}
  where $B$ is a second order linear operator defined by
\begin{align*}
B(\tilde u):= &A^{00}\partial_{00}\tilde u +\sum _j A^{0j}\partial_0\partial_j\tilde u+\sum_{i, j}\left(-\frac 13 \sum_{k, l} R_{ikjl}x_k x_l +A^{ij}\right)\partial_i\partial_j\tilde u\nonumber \\ &+B^0\partial_0 \tilde u+\sum_j\left(\sum_k\left(\frac 23 R_{ijik}+R_{0j0k}\right)x_k+B^j\right)\partial_j\tilde u
\end{align*}
where the Riemann tensor $R_{ijkl}$ and the metric $g$ are computed along $\Gamma$, depending only on $x_0$, while the function $A^{\alpha \beta}$ and $B^{\alpha}$ do depend on $(x_0, x)$ and enjoy the following decomposition:
$$A^{00}=\sum_{k, l} A^{00}_{kl}x_k x_l; \qquad A^{ij}=\sum_{k, l, m}A^{ij}_{kl}x_k x_l x_m;\qquad A^{0j}=\sum_{k, l}A^{0j}_{kl}x_k x_l$$ $$B^0=\sum_k B^0_k x_k;\qquad B^j=\sum_{k, l}B^j_{kl}x_k x_l$$ where $A^{00}_{kl}$, $A^{ij}_{kl}$, $A^{0j}_{kl}$, $B^0_k$ and $B^j_{kl}$ are smooth functions depending on $(x_0, x)$.
\end{lemma}
\begin{proof}
We argue exactly as in Section 4 of \cite{DMP} taking into account the following expansion of the metric $g$ in a neighborhood of the geodesic
\begin{equation}\label{metrica}\left\{\begin{aligned}
&g_{00}(x)=1+\sum_{k, l=1}^N R_{0k0l}x_k x_l + O(|x|^3)\\
&g_{0j}(x)=O(|x|^2), \quad  j=1, \ldots, N.\\
&g_{ij}(x)=\delta_{ij}+\frac 12 \sum_{k, l} R_{ikjl} x_k x_l + O(|x|^3) ,\quad i, j=1, \ldots, N.
\end{aligned}\right.\end{equation} whose proof is postponed in the Appendix.\end{proof}

 \subsection{The scaled problem} We write an approximated solution of problem \eqref{3pal}.
 Let
\begin{equation}\label{ansa}
\tilde u_\epsilon(x_0,x)=\mu_\epsilon(x_0)^{-\frac{N-2}{2}}w\left( {x-d_\epsilon(x_0)\over\mu_\epsilon(x_0)} \right) ,\end{equation}
where the bubble $w$ is defined in \eqref{stanbu}, {
and $d_\epsilon$ satisfies
\begin{equation}
\label{periodicA}
d_\epsilon(0) = A d_\epsilon(2\ell),\ \hbox{with}\ d_\epsilon(x_0)=\(d_{\epsilon 1}(x_0),\dots,d_{\epsilon N}(x_0)\)
\end{equation}
and $A=(a_{ij})$ is the matrix defined by \eqref{defA}.
In the sequel, $C^2_{2\ell}(\mathbb R,\mathbb R^N)$ is the space of functions   $d:[0,2\ell]\to \mathbb R^N$ which satisfy \eqref{periodicA}.}

We will take  $d_\e(x_0)$ of the form
\begin{equation}\label{di}
d_{\e j}(x_0)=\e d_j(x_0)\ \hbox{with}\ d_j\in C^2_{2\ell}(\mathbb R),\ j=1,\dots,N\
\end{equation}
and the concentration parameter $\mu_\e(x_0)$ is given by
\begin{equation}\label{mu}
\mu_\e(x_0)=\sqrt{\e} \tilde{\mu}_\e(x_0),\ \tilde{\mu}_\e(x_0)=\mu_0(x_0)+(\e\ln\e)\mu_1(x_0)+\e \mu(x_0) ,
\end{equation}
with $ \mu_0,\mu_1,\mu\in C^2_{2\ell}(\mathbb R)$.
We point out that
in \eqref{mu} and \eqref{di} the  $\mu_0$, $\mu_1$, $\mu$ and $d_j$, $j=1, \ldots, N$ are unknown  functions  which will be found    in the final step of the infinite-dimensional reduction. In particular, it will turn out that $\mu_0$ is a non degenerate solution to problem \eqref{ip-meno} in the subcritical case
or to problem \eqref{ip-piu} in the supercritical case. \\

 Therefore, it is natural to consider the change of variables
\beq\label{change}
\tilde{ u_\epsilon} (x_0, x)=\mu_\e^{-\frac{N-2}{2}}v\left(\frac{x_0}{\rho}, \frac{x-d_\e}{\mu_\e}\right),\ \rho:=\sqrt{\e},.
\eeq
Here $v_\epsilon=v_\epsilon(y_0, y)$ is defined in a region of the form
\begin{equation}\label{ddeta}
\D =\left\{(y_0, y)\ :\  y_0\in\[-{\ell\over\rho},{\ell\over\rho}\],\quad |y|<\frac{\eta}{\sqrt{\rho}}\right\}.\end{equation}

It is clear that if $\tilde u_\e(x_0,x)$ solves equation \eqref{3pal}, then    $v_\epsilon=v_\epsilon(y_0, y)$ solves problem
\beq\label{pA}
\left\{
\begin{array}{lr}
  {\mathcal A}(v)-\mu_\e^2 h v +\mu_\e^{\pm\frac{N-2}{2}\e}f_\e(v)=0\ \mbox{in}\  \D  \\
v\left(y_0+\frac{2\ell}{\rho}, y\right)=
v(y_0, Ay)
\  \mbox{for any}\  (y_0, y)\in \D.
\end{array}
\right.
\eeq
We agree that we take
$\mu_\e^{+\frac{N-2}{2}\e}$ in the supercritical case, i.e.   $f_\e(s)=(s^+)^{p+\e}$ and $\mu_\e^{-\frac{N-2}{2}\e}$  in the subcritical case, i.e. $f_\e(s)=(s^+)^{p-\e}.$\\
In \eqref{pA}
$ {\mathcal A}$ is a second order operator of the form defined in the following Lemma, whose proof can be obtained arguing exactly as in Lemma 5.1 of \cite{DMP}.

 \begin{lemma}\label{lemmalaplace}
After the change of variable \eqref{change}, the following holds true:
$$ \A(v):=a_0 \partial_{00}v+\Delta_y v +\tilde \A(v), $$ with
\begin{equation}\label{a0}
a_0(\rho y_0)=\rho^{-2}\mu_\e (\rho y_0)^2=(\mu_0+\rho\mu)^2 \end{equation}
and
$\tilde \A(v):=\sum_{\kappa=0}^2 \A_\kappa(v)+{\mathcal B}(v)$
where
\begin{align*}
\A_0(v) =& \dot{\mu}_\e^2\left[D_{yy}v[y]^2+N D_y v[y]+\frac{N(N-2)}{4}v\right] +\dot{\mu}_\e\left[D_{yy}v[y]+\frac{N-2}{2}D_y v\right][\dot{d}_\e]\\ &+D_{yy}v[\dot{d}_\e]^2 -2\mu_\e\left[\rho^{-1}D_y (\partial_0 v)[\dot{\mu}_\e y+ \dot{d}_\e]+\frac{N-2}{2}\dot{\mu}_\e\rho^{-1}\partial_0 v\right]-\mu_\e D_y v [\ddot{d}_\e]\\ &-\mu_\e\ddot{\mu}_\e\left(\frac{N-2}{2}v+ D_y v[y]\right)
\\ \A_1(v):=&-\frac 13 \sum  R_{ikjl}\left(\mu_\e y_k +d_{\e k}\right)\left(\mu_\e y_l + d_{\e l}\right)\partial_{ij} v\\   \A_2(v):=&\sum \left(\frac 23 R_{ijik}+R_{0j0k}\right)\left(\mu_\e y_k + d_{\e k}\right)\mu_\e \partial_j v  \end{align*}
and the operator $B(v)$ satisfies
\begin{align*}
{\mathcal B}(v) =& O\left(|\mu_\e y+d_\e|^2\right)\A_0(v)+ O\left(|\mu_\e y+d_\e|^3\right)\partial_{ij}v\\
 &+O\left(|\mu_\e y +d_\e|^2\right)  \left[ \mu_\e\rho^{-1}\partial_{0j}v+\mu_\e \rho^{-1}\partial_0 v - D_y (\partial_j v)[d_\e]\right.\\
 &\hskip3truecm  -\left(\frac{N-2}{2}\partial_j v + D_y (\partial_j v)[y]\right)\dot{\mu}_\e-D_y v[\dot{d}_\e]\\ &\hskip3truecm \left. -\dot{\mu}_\e\left(\frac{N-2}{2}v+ D_y v[y]\right)+\mu_\e \partial_j v\right].
\end{align*}
\end{lemma}

Our approximation close to the geodesic is
\beq\label{app}
{\boldsymbol{\tilde\omega}}=\omega+\omega_1.
\eeq
The first order approximation $\omega$ is given in \eqref{omega}, while the second order approximation $\omega_1$ is given in \eqref{correzione}.
We also set
\beq\label{Errore}
\S_\e(v):=\A(v)-\mu_\e^2 h v +\mu_\e^{\pm\frac{N-2}{2}\e}f_\e(v).
\eeq

\subsection{The ansatz: the first order approximation}\label{primo}
 We define $\omega$ to be
\beq\label{omega}
\omega:=(1+\alpha_\e)w+e_\e(\rho y_0) \chi_\e(y) Z_0(y).
\eeq
In the first term of \eqref{omega},     $w$ is the bubble defined in \eqref{stanbu} and
$\alpha_\e:=\mu_\e^{\frac{(N-2)^2}{8}\e}-1$ in the subcritical case or $\alpha_\e:=\mu_\e^{-\frac{(N-2)^2}{8}\e}-1$ in the subcritical case.
In the second term of \eqref{omega},  $\chi_\e(y):=\chi\(\epsilon^{1\over  2}|y|\)$    where $\chi$ is a cut-off function such that $\chi(s)=1$ if $s\le \delta$ and $\chi(s)=0$ if $s\ge2\delta$  with $\delta>0$ small but fixed. Moreover, $Z_0$ denotes the first eigenfunction in $L^2(\mathbb R^N)$ of the problem (see Section \ref{lineare})
\begin{equation}\label{zo}
\Delta Z_0+pw^{p-1}Z_0=\lambda_1Z_0\ \mbox{in}\ \mathbb R^N,\quad \hbox{with}\ \lambda_1>0\ \hbox{and}\ \int_{\mathbb R^N}Z_0^2\, dy=1.\end{equation}
Finally,   the function $e_\e(x_0)$ is given by
\beq\label{ee}
e_\e= \e \te_\e,\quad  \te_\e= e_0+ (\e\ln\e)e_1+\e e,
\eeq
 with $e_0,e_1,e\in C^2_{2\ell}(\mathbb R).$ We point out that $e_0,$ $e_1$ and $e$ are unknown  functions  which   will be chosen   in the final step of the infinite-dimensional reduction, together with the functions $\mu_0,$ $\mu$ and  $d_j$ introduced in \eqref{di} and \eqref{mu}. \\

Let us  estimate  the error $ {S_\e(\omega)}$  one commits by considering $\omega$ a real solution to \eqref{pA}, which is itself a function of the parameter functions $\mu, d, e$.\\ Assume that the functions $\mu, d,e $  defined respectively in \eqref{mu}, \eqref{di} and \eqref{ee},  satisfy the assumption
\beq\label{norma}
\|(\mu, d, e)\|:=\|\mu\| + \|d \| +\|e\|_\e\leq C
\eeq
for some constant $C>0$, independent of $\e$, where
\begin{align}\label{normamu}
&\|\mu\| :=\|\ddot{\mu}\|_\infty+\|\dot{\mu}\|_\infty+\|\mu\|_\infty,\ \|d \| :=\sum\limits_{j=1}^N\|d_j\|_\infty, \\ \label{normae}
&\|e\|_\epsilon:=\|\e\ddot{e}\|_\infty+\|\e^{\frac 12}\dot{e}\|_\infty+\|e\|_\infty
\end{align}
Here and in the rest of the paper, the dot denotes the derivative with respect to $x_0$.\\
It is possible to compute the expansion of the error $ {S_\e(\omega)}$ as showed in the following Lemma whose proof is postponed in Section \ref{app-a}.
\begin{lemma}\label{errore1}
If $\e>0$ small enough, then for any $(y_0, y)\in \D$ the following expansion holds
\begin{align}\label{experr}
\nonumber
 &\S_\e(\omega)= \pm\e w^p \ln w+\e \lambda_1 e_0 Z_0-\e\mu_0^2 h w +\\
\nonumber
& +\e\left[\dot{\mu}_0^2 \left(D_{yy}w[y]^2+ N D_y w[y]+ \frac{N(N-2)}{4}w\right)-\mu_0\ddot{\mu}_0 Z_{N+1}+\right.\\
\nonumber
& \left.+\mu_0^2\left(-\frac 13 R_{ikjl}y_k y_l \partial_{ij}w+\left(\frac 23 R_{ijik}+R_{0j0k}\right)y_k \partial_j w\right)\right]\\
\nonumber
& +\e^{\frac 32}\left[-\mu_0 \partial_j w \ddot{d}_j-\frac 13 \mu_0 R_{ikjl}y_k y_l \partial_{ij}w+\mu_0\left(\frac 23 R_{ijik}+R_{0j0k}\right)d_k \partial_j w -2\dot{\mu}_0 \partial_j Z_{N+1}\dot{d}_j\right]\\
\nonumber
& +\e^2\left[\left(\rho^2 a_0\ddot{e}+\lambda_1 e\right)Z_0+\left(\sum_{i, j}\dot{d}_i\dot{d}_j -\frac 13 R_{ijkl} d_k d_l \right)\partial_{ij}w+\Upsilon_0\right.+\\
\nonumber
& \left.-2\mu_0 \mu h w+ b(\rho y_0, \mu, d, e) w^p+2\dot{\mu}_0\dot{\mu}\left(D_{yy}w[y]^2+ N D_y w[y]+ \frac{N(N-2)}{4}w\right)+\right.\\
\nonumber
& \left.-\mu_0\ddot{\mu}Z_{N+1}-\mu\ddot{\mu}_0 Z_{N+1}+2\mu_0 \mu \left(-\frac 13 R_{ikjl}y_k y_l \partial_{ij}w+\left(\frac 23 R_{ijik}+R_{0j0k}\right)y_k \partial_j w\right)\right.+\\
\nonumber
& \left. -e_0 \ddot{\mu}_0 \mu_0 Z_{N+1}+\mu_0^2 e_0 \left(-\frac 13 R_{ikjl}y_k y_l \partial_{ij}Z_0+\left(\frac 23 R_{ijik}+R_{0j0k}\right)y_k \partial_j Z_0\right)+\right.\\
\nonumber
& \left.+\dot{\mu}_0^2\left(D_{yy}Z_0[y]^2+ N D_y Z_0[y]+\frac{N(N-2)}{4}Z_0\right)-\mu_0^2 h Z_0\right]\\
\nonumber
& +\e^{\frac 52}\left[-\mu \partial_j \ddot{d}_j-\frac 13 \mu R_{ikjl}y_k d_l \partial_{ij}w-\mu\left(\frac 23 R_{ijik}+R_{0j0k}\right) d_k \mu \partial_j w-2\dot{\mu}\partial_j Z_{N+1}\dot{d}_j\right.\\
\nonumber
& \left. -\mu_0 e_0 \partial_j Z_0 \ddot{d}_j-\frac 13 \mu_0 e_0 R_{ikjl}y_k d_l \partial_{ij}Z_0+\mu_0 e_0\left(\frac 23 R_{ijik}+R_{0j0k}\right) d_k \partial_j Z_0+\right.\\
& \left.-2\dot{\mu}_0 e_0\left(\frac{N-2}{2}D_y Z_0+D_{yy}Z_0[y]\right)[\dot{d}]\right]+\e^3 \Theta
\end{align}
where
\begin{itemize}
\item[-] $Z_0$ is defined in \eqref{zo} and $Z_{N+1}$ is defined in \eqref{zn}
\item[-]  the first term is   "$-\e w^p \ln w$"  in the subcritical case or "$+\e w^p \ln w$" in the supercritical case.
\item[-]
  \begin{equation}\label{ups0}
 \Upsilon_0=p(p-1)e_0^2 w^{p-2}Z_0^2+pe_0 w^{p-1}\ln w Z_0\end{equation}
 \item[-] $\Theta=\Theta(y_0,y)$ is a sum of functions of the form $$h_0(\rho y_0) \left[f_1(\mu, d,   \dot{\mu}, \dot{d})+o(1) f_2(\mu, d, e, \dot{\mu}, \dot{d}, \dot{e}, \ddot{\mu}, \ddot{d}, \ddot{e})\right]f_3(y)$$ with
\begin{itemize}
\item[-] $h_0$ a smooth function uniformly bounded in $\e$
\item[-] $f_1$ and $f_2$   smooth functions of their arguments, uniformly bounded in $\e$ when $\mu, d$ and $e$ satisfy \eqref{norma}
\item[-] $f_2$ depending linearly on the argument $(\ddot{\mu}, \ddot{d}, \ddot{e})$
\item[-] $o(1)\to0$ as $\epsilon\to0$ uniformly when $\mu, d$ and $e$ satisfy \eqref{norma}
\item[-]  $\sup_{y\in \mathbb R} (1+|y|^{N-2}) |f_3(y)|<+\infty $
\end{itemize}\end{itemize}
\end{lemma}

Now, we   use formula \eqref{experr} to compute, for each $y_0\in[-\ell/\rho,+\ell/\rho]$, the $L^2(\D_{y_0})$ the projection of the error $\S_\e(\omega)$ along the elements of the kernel of the linear operator  $\mathcal L_0:=\Delta_{\mathbb R^N}+pw^{p-1} I$  (see Section \ref{lineare}), i.e. the functions
\begin{equation}\label{zn}
Z_k(y):=\partial_{k}w(y),\  k=1, \ldots, N\ \hbox{and}\ Z_{N+1}(y):=y\cdot\nabla w(y)+\frac{N-2}{2}w(y).\end{equation}
\\
\begin{lemma}\label{errore2}
If $\e>0$ small enough, then for any $x_0=\rho y_0$ with $y_0\in[-\ell/\rho,+\ell/\rho]$ the following expansion hold:
\begin{align*}
\int_{\D_{y_0}}\S_\e(\omega)Z_k\, dy=&\e^{\frac{3}{2}}c_1 \mu_0\left(-\ddot{d}_k+\sum R_{0k0l}d_l\right) +\e^2 \theta,\quad\hbox{for any $k=1, \ldots, N$};
\end{align*}
 moreover,  if $\mu_0$ solves either \eqref{ip-meno} or \eqref{ip-piu}
  there exist  $ \mu_1,e_0,e_1\in C^2_{2\ell}(\mathbb R)$   such that
\begin{align*}
\int_{\D_{y_0}}\S_\e(\omega)Z_{N+1}\, dy=&\e^2 c_2\mu_0\[\alpha_{N+1}(x_0)+ c_3Q( x_0,d)   - \ddot \mu  +\(a_n\sigma \mp {b_n\over\mu_0^2}\)\mu\]   +\e^3 |\ln\e|\theta
\end{align*}
 and
\begin{align*}
\int_{\D_{y_0}}\S_\e(\omega)Z_0\, dy =& \e^2  \left[\e a_0 \ddot{e}+\lambda_1 e +\alpha_0(x_0)+c_4 Q( \rho y_0,d) +  \beta (x_0)\mu\right] +\e^3|\ln\e| \theta.
\end{align*}
Here
\begin{itemize}
\item[-] $\sigma $ is defined in \eqref{gamma0} and $a_n,\ b_n$ are   positive constants depending only on $n$ defined in \eqref{costanti}
        \item[-] $Q( x_0,d):= \sum  \left(\dot{d}_j^2-\frac 13 R_{ikjl}d_k d_l\right)   $
       \item[-] $c_i$'s are constants which depends only on $n$
        \item[-]   $\alpha_i$'s and $\beta$ are explicit smooth functions, uniformly bounded in $\e$ when $\mu, d$ and $e$ satisfy \eqref{norma}
             \item[-]
      $\theta=\theta(x_0)$
      denotes a sum of functions of the form $$h_0(x_0)\left[h_1(\mu, d, e, \dot{\mu}, \dot{e}, \dot{d})+o(1)h_2(\mu, d, e, \dot{\mu}, \dot{d}, \dot{e},\ddot{\mu}, \ddot{d}, \ddot{e})\right],$$ where
\begin{itemize}
 \item[-] $h_0$ is a smooth function uniformly bounded in $\e$
   \item[-]$h_1$ and $h_2$ are smooth functions of their arguments, uniformly bounded in $\e$ when $\mu, d$ and $e$ satisfy \eqref{norma}
      \item[-] $h_2$ depends linearly on the argument $(\ddot{\mu}, \ddot{d}, \ddot{e})$
   \item[-] $o(1)\rightarrow 0$ as $\e\rightarrow 0$ uniformly  when $\mu, d$ and $e$ satisfy \eqref{norma}
                \end{itemize}     \end{itemize}\end{lemma}
The proof is postponed in Section \ref{app-b}.\\

In the sequel we will use the following norms, which are motivated by the linear theory presented in Section~\ref{lineare}.
For functions $\phi$, $g$ defined on a set  $\D$  as in \eqref{ddeta}, and for a fixed  $2\leq \nu <N$,  let
\begin{align*}
\|\phi\|_* & :=\sup_{\D}(1+|y|^{\nu-2})|\phi(y_0, y)|+\sup_{\D}(1+|x|^{\nu-1})|D\phi(x_0, x)| ,
\\
\|g\|_{**} & :=\sup_{\D}(1+|y|^\nu)|g(y_0, y)|.
\end{align*}

Therefore,     from the expansion given in \eqref{experr} we conclude that the error $\S_\e(\omega)$, computed in \eqref{experr}, has the properties listed in the following Lemma.

   \begin{lemma}
 Let $\mu_0$ and $e_0$  as in Lemma \ref{errore2}
   If $\e$ is small enough
\beq\label{experr1}
\S_\e(\omega)=\e S_0 +\e\left[\rho^2 a_0 \ddot{e}+\lambda_1 e \right] \chi_\e Z_0 +N_0
\eeq
where
\begin{itemize}
\item[-] $S_0$ is a smooth function of $\rho y_0$  uniformly bounded in $\e$
\item[-] $S_0$ does not depend on $\mu, d$ and $e$.
 \item[-] $\int_{\D_{y_0}}S_0 Z_j\, dy =0$ for any $y_0\in (-\rho^{-1}\ell,\rho^{-1}\ell)$  and
 for any $j=0, \ldots, N+1$
 \item[-] $\|N_0\|_{**}\leq c\e^{\frac 32}$
 \end{itemize}
Here $c$ is a positive constant   independent of $\e$.
All the estimates are uniform with respect to $\mu, d$ and $e$ which satisfy \eqref{norma}.
 \end{lemma}

\subsection{The ansatz:  the second order approximation}\label{secondo}
Now we introduce a further correction $  \omega_1 $ to $\omega$, to get the final approximation ${\boldsymbol {\tilde\omega}}:=\omega+\omega_1 $. The correction $\omega_1 $ is chosen to reduce the size of the error \eqref{experr1}, killing the term $\e S_0$ and it is found in the following Lemma, whose proof  can be carried out arguing exactly as in Section 5 of \cite{DMP}.

\begin{lemma}\label{corr-omega1}
If $\e$ is small enough there exists a unique solution $\omega_1$ of the problem
\beq\label{correzione}
\left\{
\begin{array}{lr}
 {\A}(\omega_1 )-\mu_\e^2 h \omega_1 +p w^{p-1}\omega_1=-\e S_0 +\sum_{j=0}^N \sigma_j Z_j \qquad \mbox{in}\quad \D\\\\
\int_{\D_{y_0}} \omega_1 (y_0, y)Z_j  dy=0\qquad \hbox{for any}\ y_0\in\[-{\ell\over\rho},{\ell\over\rho}\],\   j=0, \ldots, N+1
\end{array}
\right.
\eeq
Moreover, the function $ \omega_1$ satisfies
\begin{itemize}
\item[-] $\|\omega_1\|_*\leq c\e$ and $\|\partial_0\omega_1\|_*\leq c\e^{\frac 32}$
\item[-] $\omega_1$ depends smoothly on   $\mu$ and $d$ and it is independent on $e$
\item[-] $ \|\omega_1(\mu_1, d_1)-\omega_1(\mu_2, d_2)\|_*\leq c \|(\mu_1-\mu_2, d_1-d_2)\|$
\end{itemize}
and each    function $\sigma_j$ satisfies
\begin{itemize}
\item[-] $ \|\sigma_j\|_\infty\leq o(1)\e^3$
\item[-]   $\sigma_j$ depends smoothly on   $\mu$ and $d$ and it is independent on $e$
\item[-] $ \|\sigma_j(\mu_1, d_1)-\sigma_j(\mu_2, d_2)\|_\infty\leq c\e^2 \|(\mu_1-\mu_2, d_1-d_2)\|$
\end{itemize}
Moreover, it holds true
\beq\label{experr2}
\S_\e({\boldsymbol {\tilde\omega}})=\e^{\frac 32} S_1 +\e\left[\rho^2 a_0 \ddot{e}+\lambda_1 e \right] \chi_\e Z_0 +N_1+\sum_{j=0}^N \sigma_j Z_j
\eeq
where
\begin{itemize}
\item[-] $S_1$ is a smooth function of $\rho y_0$  uniformly bounded in $\e$
\item[-] $S_1$   depends smoothly on $\mu, d$ and $e$.
\item[-] $ \|S_1 (\mu_1, d_1,e_1)-S_1(\mu_2, d_2,e_2)\|_{**}\leq c  \|(\mu_1-\mu_2, d_1-d_2,e_1-e_2)\|$
 \item[-] $\|N_1\|_{**}\leq c\e ^2$
 \end{itemize}
 Here $c$ is   positive constant independent of $\e$. All the estimates are uniform with respect to $\mu, d$ and $e$ which satisfy \eqref{norma}.
 Moreover, the components of $\S_\e({\boldsymbol {\tilde\omega}})$ along the $Z_j$'s satisfy the estimate in Lemma \ref{errore2}.
 \end{lemma}

\section{The  error $ {S_\e(\omega)}$}\label{sec-errore}
\subsection{The pointwise estimate of the error  }\label{app-a}
We recall that
$$ \S_\e(\omega)= \A(\omega)-\mu_\e^{2} h \omega +\mu_\e^{\pm\frac{N-2}{2}\e}f_\e(\omega)  $$
where by Lemma \ref{lemmalaplace}
$$\A(\omega)=a_0\partial_{00}\omega+\Delta_y \omega+\underbrace{\sum_{k=0}^2\A_k(\omega)+{\mathcal B}(\omega)}_{\tilde{\A}(\omega)}$$ and $$\omega(y)=(1+\alpha_\e)w(y)+e_\e(\rho y_0) \chi_\e (y) Z_0(y).$$ Here we recall that
$$ \alpha_\e=\mu_\e^{\mp\frac{(N-2)^2}{8}\e}-1  $$ and
$$\Delta\((1+\alpha_\e)w\)+\mu_\e^{\pm\frac{N-2}{2}\e}f_0\((1+\alpha_\e)w\)=0\qquad \mbox{in}\,\, \mathbb R^N .$$

\begin{proof}[Proof of Lemma \ref{errore1}]
We use Lemma \ref{lemmalaplace}.\\ A straightforward computation shows that
\begin{align}
\S_\e(\omega)&=  \underbrace{\sum_{\kappa=0}^2 \A_\kappa(w) -\mu_\e^2    h w\pm\e w^p \ln w+\left[\rho^2 a_0 \ddot{e}_\e(\rho y_0)+\lambda_1 e_\e(\rho y_0)\right]\chi_\e  Z_0}_{J_0}\nonumber \\ & +\underbrace{{\mathcal B}(w)+  a_0w\partial_{00}\alpha_\e +\tilde\A(\alpha_\e w)-\mu_\e^2  \alpha_\e  h w}_{J_1}\nonumber\\ &
+ \underbrace{\mu_\e^{\pm\frac{N-2}{2}\e}\left[ f_\e\left((1+\alpha_\e)w\right) -f_0\left((1+\alpha_\e)w\right) \right] \mp\e w^p \ln w}_{J_2}\nonumber\\ &
+\underbrace{\sum_{\kappa=0}^2 \A_\kappa(e_\e \chi_\e Z_0) -\mu_\e^2    e_\e \chi_\e Z_0h}_{J_3}
\nonumber \\ &  +\underbrace{ {\mathcal B}(e_\e \chi_\e Z_0)+  e_\e  Z_0 \Delta \chi_\e +2 e_\e \nabla \chi_\e\nabla Z_0
  }_{J_4} \nonumber  \\
 &  + \underbrace{\mu_\e^{\pm\frac{N-2}{2}\e}\left[f_\e(\omega )-f_\e\left((1+\alpha_\e)w\right)\right]-f' _0(w)e_\e \chi_\e Z_0}_{J_5}.
\end{align}

By Lemma \ref{lemmalaplace}, we get the first term of $J_0$
\begin{align}\label{I2}
\sum_{\kappa=0}^2 \A_\kappa(w)&= \dot{\mu}_\e^2\left[D_{yy}w[y]^2+N D_y w[y]+\frac{N(N-2)}{4}w\right]\nonumber\\
 &+\dot{\mu}_\e\left[D_{yy}w[y]+\frac{N-2}{2}D_y w\right][\dot{d}_\e]+D_{yy}w[\dot{d}_\e]^2\nonumber\\
 &-\mu_\e D_y w [\ddot{d}_\e]-\mu_\e\ddot{\mu}_\e\left(\frac{N-2}{2}w+D_y w[y]\right)\nonumber\\  & -\frac 13 \sum R_{ikjl}\left(\mu_\e y_k +d_{\e k}\right)\left(\mu_\e y_l + d_{\e l}\right)\partial_{ij} w\nonumber\\
 &+\sum\left(\frac 23 R_{ijik}+R_{0j0k}\right)\left(\mu_\e y_k + d_{\e k}\right)\mu_\e \partial_j w+\e^3\Theta\nonumber\\
&=  \e^2\left[\sum \left(\dot{d}_i\dot{d}_j-\frac 13 R_{ikjl}d_k d_l\right) \right]\partial_{ij}w \nonumber\\ &  +\rho \e \left[-\tm D_y w [\ddot{d}]-\sum \frac 13 \tm R_{ikjl}y_k d_l \partial_{ij}w+\right.\nonumber\\
& \left. \qquad+\left(\frac 23 R_{ijik}+R_{0j0k}\right)d_k \tilde{\mu}\partial_j w-2\dot{\tilde{\mu}}D_y Z_{N+1}[\dot{d}]\right]\nonumber\\
& +\rho^2\left[\dot{\tilde{\mu}}^2\left[D_{yy}w[y]^2+N D_y w[y]+\frac{N(N-2)}{4}w\right]-\tilde{\mu}\ddot{\tilde{\mu}}Z_{N+1}\right.\nonumber\\
& \left.\qquad+\tilde{\mu}^2\left(-\frac 13 \sum R_{ikjl} y_k y_l \partial_{ij}w+\left(\frac 23 R_{ijik}+R_{0j0k}\right) y_k \partial_j w\right)\right]+\e^3\Theta,
\end{align}
where $\Theta=\Theta(\rho y_0,y)$ has the required properties.\\
By Lemma \ref{lemmalaplace}, we  deduce that $\mathcal B(w)$ is of lower order with respect to $\sum \A_k(w)$.  Moreover, by definition of $\alpha_\e$ we get that $\alpha_\e=O(\e|\ln\e|)$ as $\e\rightarrow 0$. Hence $\alpha_\e \tilde{\A}(w)$ and $\mu_\e\alpha_\e h w$ are terms of lower order with respect to the others. Furthermore $\partial_{00}\alpha_\e=\rho^2 O(\alpha_\e)$, so also $a_0 \partial_{00}[\alpha_\e w]=O(\e^2|\ln\e|)w$.
Therefore,
 $$J_1= \e^3 \Theta$$
 where  $\Theta=\Theta(\rho y_0,y)$  is a sum of functions of the form $h_0(\rho y_0) f_1(\mu, d, \dot{\mu}, \dot{d})f_2(y)$, with $h_0$ a smooth function uniformly bounded in $\e$, $f_1$ a smooth function of its arguments, homogeneous of degree $3$, uniformly bounded in $\e$ and $ \sup_{y\in\mathbb R} (1+|y|^{N-2})|f_2(y)|<+\infty.$

By mean value theorem we deduce that
\begin{align}\label{jei2}
&J_2=\pm{(n-2)^2\over 8}(\e^2\ln\e )w^p(\ln w-1)\pm\e^2w^p\({(n-2)^2\over 8} ( \ln w-1)\ln\mu+{1\over2} \ln w\)\nonumber \\ &+O\(\e^3|\ln\e|\).
 \end{align}

By Lemma \ref{lemmalaplace} we  also get   that
\begin{align*}
J_3&= \e \te \left\{\e^2\left[\left(\sum \dot{d}_i\dot{d}_j -\frac 13 R_{ikjl}d_k d_l\right)\partial_{ij}Z_0\right]\right.\\ &+\rho\e\left[-\tm D_y Z_0[\ddot{d}]-\frac 13 \tm R_{ikjl}y_k d_l \partial_{ij}Z_0 +\tm\left(\frac 23 R_{ijik}+R_{0j0k}\right)d_k \partial_j Z_0\right.\\
&\qquad\left.-2\dot{\tm}\left(\frac{N-2}{2}D_y Z_0+D_{yy}Z_0[y]\right)[\dot{d}]\right]\\
& \quad+\rho^2\left[-\ddot{\tm}\tm Z_{N+1}+\tm^2\left(-\frac 13 R_{ikjl}y_k y_l \partial_{ij}Z_0+\left(\frac 23 R_{ijik}+R_{0j0k}\right)y_k\partial_j Z_0\right)+\right.\\
& \quad\left.\left.+\dot{\tm}^2\left(D_{yy}Z_0[y]^2+ N D_y Z_0[y]+\frac{N(N-2)}{4}Z_0\right)-\tm^2 h Z_0\right]\right\}
\\
&+\rho\e \dot\te\left\{\e\left(-2\tm D_y Z_0[\dot{d}]\right)+\rho\e\left[-2\tm\dot{\tm}D_y Z_0[y]-(N-2)\tm\dot{\tm}Z_0\right]\right\} \end{align*}
and
$$J_4=\e^3 \Theta$$ where $\Theta=\Theta(\rho y_0,y)$ has the required properties.\\
Finally, standard estimates yield to
$$J_5=\e^2\underbrace{\left[p(p-1)e_0^2 w^{p-2}Z_0^2+p e_0 w^{p-1}\ln w Z_0\right]}_{\Upsilon_0}+\e^3|\ln\e|\Theta,
$$
where  $\Theta=\Theta(\rho y_0,y)$  is a sum of functions of the form $h_0(\rho y_0)h_1(\mu, d, e)h_2(y)$ with $h_0$ a smooth function, uniformly bounded in $\e$, $h_1$ a smooth function of its arguments and $\sup _{y\in\mathbb R}(1+|y|^{N-2})|h_2(y)|<+\infty$.\\  Collecting all the previous estimates we get the proof.\end{proof}

\subsection{The components  of the error  along the $\boldsymbol{Z_j}$'s}\label{app-b}
\begin{proof}[Proof of Lemma \ref{errore2}]
The proof consists of two steps. In the first part we compute the expansion in $\e$ of the projection assuming that
$$\mu_\e= \rho\tilde{\mu},\qquad d_{\e j}=\e d_{j},\qquad e_\e=\e \tilde{e}.$$
In the second part we will choose the $\e-$order terms $\mu_0$ and $e_0$ and the $\e\ln\e$-order terms $\mu_1$ and $e_1$  in the expansion of $\tm$ and $\te.$\\\\
Arguing as in the proof of Lemma \ref{errore1}, we have
 \begin{eqnarray*}
\S_\e(\omega)&=&\underbrace{\pm\e w^p\ln w-\rho^2\tm^2 h w}_{I_1}+\underbrace{\sum_{k=0}^2 \A_k(w)}_{I_2}+\underbrace{\e\left[\rho^2 a_0 \ddot{\te}+\lambda_1\te\right]\chi_\e Z_0}_{I_3}+\underbrace{J_1+\dots+J_5}_{I_4}.
\end{eqnarray*}
We stress the fact that the first term in $I_1$ is    $  "+\e  w^p\ln w"$ in the super-critical case and $"-\e w^p\ln w"$ in the sub-critical case.
\begin{itemize}
\item {\em The projection of $I_1$.}\\
 \begin{align*}
 \int_{\D_{y_0}}I_1Z_{N+1}\, dy&=\pm\e \int_{\D_{y_0}}w^p \ln w Z_{N+1}\, dy-\rho^2\tm^2  \int_{\D_{y_0}} h w Z_{N+1}\, dy\\
&=-\e A_1 +O(\e \rho^N)-\rho^2\tm^2h(\rho y_0)\int_{\mathbb R^N} w Z_{N+1}\, dy+ O(\rho^N)\\
& = \e\[\pm A_1 - \tm^2 h(\rho y_0) A_2\]+ O(\rho^N).\end{align*}
  where
\begin{equation}\label{a1}A_1=\int_{\mathbb R^N}w^p \ln w Z_{N+1}\, dy = \frac{N}{(p+1)^2}\int_{\mathbb R^N}w^{p+1}\, dy>0\ \hbox{(see Remark \ref{gamma})}\end{equation}
and
\begin{equation}\label{a2}A_2=\int_{\mathbb R^N}w  Z_{N+1}\, dy <0 \ \hbox{(see Remark \ref{gamma})}.\end{equation}\\
\begin{align*} \int_{\D_{y_0}} I_1  Z_{k}\, dy&= \e\int_{\D_{y_0}}w^p \ln w Z_j\, dy +\rho^2\tm^2  \int_{\D_{y_0}}hw Z_j\, dy\\ &=\e\int_{\mathbb R^N} w^p \ln w Z_j\, dy +\rho^2\tm^2  h(\rho y_0)\int_{\mathbb R^N}w Z_j\, dy+ O(\rho^{N+1})\\ &=O(\rho^{N+1})\ \hbox{for}\ k=1,\dots,N.\end{align*}\\
 \begin{align*}
\int_{\D_{y_0}}I_1 Z_0\, dy&= -\e\int_{\D_{y_0}} w^p\ln w Z_0\, dy-\rho^2\tm^2  \int_{\D_{y_0}}hwZ_0\, dy\nonumber \\ &=  \e\[- A_3- \tm^2 h(\rho y_0) A_4\]+O(\rho^N),
\end{align*}
where \begin{equation}\label{a34}A_3:=\int_{\mathbb R^N} w^p \ln w Z_0\, dy,\ A_4:=\int_{\mathbb R^N}w Z_0\, dy.\end{equation}
\item {\em The projection of $I_2$.}\\
We use estimate \eqref{I2}.
\begin{align*}\label{prg1zn1}
\nonumber
\int_{\D_{y_0}} I_2 Z_{N+1}\, dy&= \e^2\sum \left(\dot{d}_i\dot{d}_j-\frac 13 R_{ikjl}d_k d_l \right)\int_{\D_{y_0}} \partial_{ij} w Z_{N+1}\, dy\\
\nonumber
& -\rho\e \tm \sum  \ddot{d}_j \int_{\D_{y_0}} \partial_j w Z_{N+1}\, dy \\
\nonumber
&  -\frac 13\tm \rho\e\sum R_{ikjl}d_l \int_{\D_{y_0}} y_k \partial_{ij}w Z_{N+1}\\
\nonumber
& +\rho \e \tm \sum \left(\frac 23 R_{ijik}+R_{0j0k}\right) d_k \int_{\D_{y_0}} \partial_j w Z_{N+1}\, dy\\
\nonumber
&   -2\dot{\tm}\rho \e \sum  \dot{d}_j \int_{\D_{y_0}}\partial_j Z_{N+1}Z_{N+1}\, dy\\
\nonumber
& +\dot{\tm}^2\rho^2\int_{\D_{y_0}}\left[D_{yy}w[y]^2+N D_y w[y]+\frac{N(N-2)}{4}w\right]Z_{N+1}\, dy\\
\nonumber
&  -\tm\ddot{\tm}\rho^2 \int_{\D_{y_0}} Z_{N+1}^2\, dy\\
\nonumber
& -\rho^2 \tm^2\frac 13 \sum R_{ikjl}\int_{\D_{y_0}} y_k y_l \partial_{ij}w Z_{N+1}\, dy\\
\nonumber
&  +\tm^2 \rho^2\sum \left(\frac 23 R_{ijik}+R_{0j0k}\right)\int_{\D_{y_0}} y_k \partial_j w Z_{N+1}\, dy\\
\nonumber
&=  \e^2 \sum  \left[\dot{d}_i^2-\frac 13 R_{ikil}d_k d_l\right]\int_{\mathbb R^N} \partial_{ii}w Z_{N+1}\, dy\\
\nonumber
&  +\tm^2 \rho^2 \sum \left(\frac 23 R_{ijij}+R_{0j0j}\right)\int_{\mathbb R^N}y_j\partial_j w Z_{N+1}\, dy+\\
&   -\tm \ddot{\tm}\rho^2  \int_{\D_{y_0}} Z_{N+1}^2\\
\nonumber
&  -\frac 13 \rho^2 \tm^2\sum R_{ikjl}\int_{\mathbb R^N}y_k y_l \partial_{ij}w Z_{N+1}\, dy\\
\nonumber
&  +\e^3 \theta\\
&=  \e^2B_1 \underbrace{\sum  \left[\dot{d}_i^2-\frac 13 R_{ikil}d_k d_l\right]}_{Q(d,\rho y_0)}\\
\nonumber
&   +\e\left[\tm^2   \sum \left(\frac 23 R_{ijij}+R_{0j0j}\right)B_2 -\tm \ddot{\tm}  B_3\right]\\
\nonumber
&  +\e^3 \theta
\end{align*}
where the function $\theta=\theta(\rho y_0)$ has the required properties and
\begin{equation}\label{b123}B_1:=\int_{\mathbb R^N} \partial_{ii}w Z_{N+1}dy, \ B_2:=\int_{\mathbb R^N}y_j\partial_j w Z_{N+1}dy<0,\ B_3:=\int_{\mathbb R^N}   Z_{N+1}^2  dy.\end{equation}
Here we used the fact that
$$\sum R_{ikjl}\int_{\mathbb R^N}y_k y_l \partial_{ij}w Z_{N+1}\, dy=0,$$
because $R_{ikjl}$ is antisymmetric (i.e. $R_{ikjl}=-R_{kijl}$) and $\int_{\mathbb R^N}y_k y_l \partial_{ij}w Z_{N+1}\, dy$ is symmetric.\\
\begin{align*}
\int_{\D_{y_0}}I_2 Z_k\, dy&=  \rho \e \tm \left[-  \ddot{d}_k\int_{\mathbb R^N} Z_j^2\, dy -\frac 23 R_{iljm}d_l\int_{\mathbb R^N} y_m \partial_{ij} w Z_k\, dy \right.\\ &\left.\qquad\qquad+ \left(\frac 23 R_{ijil}+R_{0j0l}\right) d_l\int_{\mathbb R^N} Z_j^2\, dy\right]\\ \nonumber &  +\rho^2\e \theta \\
&=  \e^{3\over2} \tm B_4 \left[-\ddot{d}_k+R_{0j0l}d_l\right]+\rho^2 \e \theta,
\end{align*}
where \begin{equation}\label{b4}B_4:=\int_{\mathbb R^N} Z_j^2\, dy,\ j=1,\dots,N.\end{equation}
Here we used the fact that \begin{eqnarray*}
& &-\frac 23 R_{iljm}\int y_m \partial_{ij}w Z_k\, dy \\ & &=-\frac 23 \left[R_{ilik}\int y_k \partial_{ii}w Z_k\, dy+R_{ilki}\int y_l \partial_{ik}w Z_k\, dy+R_{kljj}\int y_j\partial_{kj}w Z_k\, dy\right] \\
&& = -\frac 13 B_4 \left[ R_{ilik}-R_{ilki}   \right] =-\frac 23 B_4 R_{ilik}  .
\end{eqnarray*}\\
 \begin{align*}
\nonumber
\int_{\D_{y_0}}I_2 Z_0\, dy &=  \e^2\left[\sum \left(\dot{d}_i^2-\frac 13 R_{ikil} d_k d_l \right)\int_{\mathbb R^N} \partial_{ii}w Z_0\, dy\right]\\\nonumber & +\tm^2\rho^2\sum \left(\frac 23 R_{ijij}+R_{0j0j}\right)\int_{\mathbb R^N}y_j \partial_{j}w Z_0\, dy \\
&   -\rho^2\tm^2\frac 13\sum R_{ikjl}\int_{\mathbb R^N}y_k y_l \partial_{ij} w Z_0\, dy+\e^3 r\\
&=  \e^2 B_5\underbrace{\sum  \left[\dot{d}_i^2-\frac 13 R_{ikil}d_k d_l\right]}_{Q(d,\rho y_0)}\\ \nonumber & +\e\tm^2 B_6\sum \left(\frac 23 R_{ijij}+R_{0j0j}\right)\\
&  +\e^3 \theta,
\end{align*}
where
\begin{equation}\label{b56}B_5:=\int_{\mathbb R^N} \partial_{ii}w Z_0\, dy,\ B_6:=\int_{\mathbb R^N}y_j \partial_{j}w Z_0\, dy .\end{equation}
Here we used the fact that
$$\sum R_{ikjl}\int_{\mathbb R^N}y_k y_l \partial_{ij}w Z_{0}\, dy=0,$$
because $R_{ikjl}$ is antisymmetric (i.e. $R_{ikjl}=-R_{kijl}$) and $\int_{\mathbb R^N}y_k y_l \partial_{ij}w Z_{0}\, dy$ is symmetric.\\

\item {\em The projection of $I_3$.}\\
$$\int_{\D_{y_0}}I_3 Z_{N+1}\, dy=o(1)\e^3
\ \hbox{and}\
\int_{\D_{y_0}}I_3 Z_k\, dy =o(1)\e^3\ \hbox{for any}\ k=1, \ldots, N,
$$
because of the symmetry and of the orthogonality of $Z_0$ with $Z_{N+1}$ and $Z_j$.
\\
$$
 \int_{\D_{y_0}}I_3 Z_0\, dy= \e   \left[\rho^2 a_0 \ddot{\te}+\lambda_1 \te\right]+o(1)\e^3
$$
because $\int_{\mathbb R^N}Z_0^2\, dy =1.$

\item {\em The projection of $I_4$.}\\

$$\int_{\D_{y_0}}I_4 Z_{N+1}\, dy = \e^2\ln\e D_1+\e^2b_{ 1}(\rho y_0) +\e^3|\ln\e| \theta
$$
 $$\int_{\D_{y_0}}I_4 Z_k\, dy = \e^2 \theta\ \hbox{for any}\ k=1, \ldots, N.
$$
$$
\int_{\D_{y_0}}I_4 Z_0\, dy = \e^2\ln\e D_2+\e^2b_{2}(\rho y_0) +\e^3|\ln\e| \theta
$$
where
$$D_1:=\pm{(N-2)^2\over16} A_1,\ D_2:=\pm{(N-2)^2\over16} A_3 \ \hbox{(see \eqref{a1} and \eqref{a34})},$$ $b_1,$ $b_2$ are explicit functions and
the function $\theta=\theta(\rho y_0)$ has the required properties .\\\\

\end{itemize}

Hence, summing up the previous calculations we conclude that
\begin{align}\label{pzn1}
\int_{\D_{y_0}}\S_\e(\omega)Z_{N+1}\, dy&= \e\underbrace{\(\pm A_1-\mu_0\ddot\mu_0 B_3 +\mu_0^2g_1\)}_{  \hbox{ the choice of $\mu_0$ $\Rightarrow$ =0}} \nonumber\\
&+\e^2\ln\e\underbrace{\(-\ddot\mu_1\mu_0  B_3+\mu_1\(-\ddot\mu_0 B_3+2\mu_0 g_1\)+D_1\)}_{  \hbox{ the choice of $\mu_1$ $\Rightarrow$ =0}}  \nonumber\\
&+\e^2 \(-\ddot\mu \mu_0  B_3+\mu \(-\ddot\mu_0 B_3+2\mu_0 g_1\) +B_1Q(d,x_0)+b_1(x_0)\)  \nonumber\\
&+O(\e^3|\ln\e|)\end{align}
where (see Remark \ref{gamma})
\begin{equation}\label{gamma1}
g_1(x_0):=-A_2 h(x_0) +\sum\left(\frac 23 R_{ijij}+R_{0j0j}\right)B_2=-A_2\sigma (x_0) \end{equation}
and
\begin{align}\label{pz0}
\int_{\D_{y_0}}\S_\e(\omega)Z_{0}\, dy&=\e\underbrace{\(\lambda_1 e_0 -A_3 +\mu_0^2 g_2\)}_{  \hbox{ the choice of $e_0$ $\Rightarrow$ =0}} \nonumber\\
&+\e^2 \ln\e\underbrace{\(\lambda_1 e_1 +2\mu_0\mu_1 +D_2\)}_{  \hbox{ the choice of $e_1$ $\Rightarrow$ =0}} \nonumber\\
&+\e^2 \(\e a_0\ddot e+\lambda_1 e+a_0\ddot e_0+b_2(x_0)+2\mu_0\mu g_2+B_5Q(d,x_0)  \)  \nonumber\\
&+O(\e^3|\ln\e|)\end{align}
where
  \begin{equation}\label{gamma2}
g_2(x_0):=- A_4 h(x_0) +\sum\left(\frac 23 R_{ijij}+R_{0j0j}\right)B_6.\eeq
\\
 More precisely, $\mu_0$ solves the periodic  O.D.E.
\beq\label{eqmu0}
-\ddot{\mu}_0 B_3 + g_1\mu_0\pm\frac{A_1}{\mu_0}=0,\   \mu_0>0\ \hbox{in}\ [0,2\ell].
\eeq
which is nothing but   problem \eqref{ip-meno} or \eqref{ip-piu}
where  (see Remark \ref{gamma})
\begin{equation}\label{costanti}
a_n:=-{A_2\over B_3}>0 \ \hbox{ and}\ b_n:= {A_1\over B_3}>0 \ \hbox{(see \eqref{a1},  \eqref{a2} and \eqref{b123}).}
\end{equation}
Moreover,
\beq\label{eqe0}
e_0=\frac{A_3-\mu_0^2 g_2}{\lambda_1 }.
\eeq
Finally, $\mu_1$ solves the periodic  O.D.E.
\beq\label{eqmu1}
-\ddot\mu_1\mu_0  B_3+\mu_1\underbrace{\(-\ddot\mu_0 B_3+2\mu_0 g_1\)}_{\displaystyle{=\mu_0 g_1\mp \frac{A_1}{\mu_0^2}}}+D_1=0\ \hbox{in}\ [0,2\ell].
\eeq
We point out that $\mu_1$ does exist, because     $\mu_0$ is a non degenerate solution of \eqref{eqmu0} (see also Lemma \ref{lnp}).
Moreover,
\beq\label{eqe1}
e_1=\frac{-2\mu_0\mu_1-D_2}{\lambda_1 }.
\eeq
 That concludes the proof.
\end{proof}

\begin{remark}\label{gamma}
It holds
\begin{itemize}
\item
 $g_1(x_0)=-A_2\sigma (x_0)$
     with $A_2<0 $ (see \eqref{a2})
     \item $A_1>0$ (see \eqref{a1})
     \item  $a_n =-{A_2\over B_3}={8(N-1)\over (N-2)(N+2)}={8(n-2)\over (n-3)(n+1)}$  (see \ \eqref{a2} and \eqref{b123})
 \item  $b_n = {A_1\over B_3}={ (N-2)^2(N-4)\over 4(N+2)}={ (n-3)^2(n-5)\over 4(n+1)}$  (see \eqref{a1}  and \eqref{b123})
     \end{itemize}
   \end{remark}
\begin{proof}It is useful to point out that
$${B_2\over A_2} ={(N-2)(N-3)\over 4(N-1)}.
$$
Indeed, if we denote by
$$I^q_p:=\int_0^{+\infty}{r^q\over (1+r)^p}dr\ \hbox{if}\ p-q>1$$
and we use the properties
$$I^q_{p+1}={p-(q+1)\over p}I^q_p\ \hbox{and}\ I^{q+1}_{p+1}={q+1\over p-(q+1)}I^q_{p+1}$$
a straightforward computation shows that
$$A_1={N\over (p+1)^2} \int_{\mathbb R^N} w ^{p+1}\, dy= \alpha^2_N{(N-2)^4\over 8N }
\omega_N I_N^{N/2}>0,$$
$$A_2=\int_{\mathbb R^N} w Z_{N+1}\, dy=-\alpha^2_N{4(N-1)(N-2)\over  N(N-4)}
\omega_N I_N^{N/2}<0,$$
$$B_2=\int_{\mathbb R^N}y_j\partial_j w Z_{N+1}\, dy=-\alpha^2_N{(N-2)^2(N-3)\over 2N(N-4)}
\omega_N I_N^{N/2}<0
$$
and
$$B_3=\int_{\mathbb R^N} Z_{N+1}^2\, dy=\alpha^2_N{(N-2)^2(N+2)\over 2N(N-4)}
\omega_N I_N^{N/2}>0,
$$
where $\omega_N$ is the measure of the sphere ${\mathbb S}^{N-1}.$ Therefore, we immediately deduce the quantities $a_n$ and $b_n ,$ taking into account that $N=n-1.$\\
 Moreover,  the scalar curvature $R_g$   in normal coordinates reads  as
\beq\label{curve}R_g(x_0)=\sum\limits_{i,j=0}^N R_{ijij}(x_0),\eeq
so it is easy to check that
$$\frac 23 \sum_{i,j=1}^N  R_{ijij}(x_0)+\sum_{j=1}^NR_{0j0j}(x_0) =\frac 43 S_g(x_0)-\frac13\sum_{j=0}^NR_{0j0j}(x_0).$$
Therefore, the claim follows.
\end{proof}

 \section{The infinite dimensional reduction }\label{gluing}
\subsection{The gluing procedure }

Here we perform a gluing procedure that reduces the full problem \eqref{pb} to the scaled problem \eqref{pA} in the neighborhood of the scaled geodesic.
\\ Since the procedure is very similar to that of \cite{DMP} we briefly sketch it.\\\\
We denote by $M_\rho$ the scaled manifold ${1\over \rho}M,$   by $z$   the original variable in $M_\rho$ and by $ \xi:=\rho z$ the corresponding point in $M.$
It is clear that the function $u(x)$ is a solution to   \eqref{pb} if and only if  the function $v(z):=\rho^{\frac{N-2}{2}}u(\rho z)$ solves the problem
\beq\label{pbgluing}
\Delta_g v -\rho^2 h v + \rho^{-\frac{N-2}{2}\e}v^{p-\e}=0\qquad \mbox{in}\,\, \M_\rho
\eeq
The function ${\boldsymbol {\tilde\omega}}(y_0, y)$ constructed in   \eqref{app} defines an approximation  to a solution of \eqref{pb} near the geodesic through the natural change of variables \eqref{change}. \\
It is useful to introduce the following   notation.
Let  $f(z)$ be a function defined in a small neighborhood of the scaled geodesic $\Gamma_\rho:=  {1\over \rho}\Gamma .$  Through the   change of variables \eqref{change} we denote by \begin{equation}\label{changeglu}\tilde{f}(y_0, y)=\tm_\e^{-\frac{N-2}{2}}(\rho y_0)f\({1\over\rho}F\(\rho y_0, \mu_\e(\rho y_0)+d_\e(\rho y_0)\)\),\end{equation}
 where the point $ \rho z=F(\rho y_0, \mu_\e(\rho y_0)+d_\e(\rho y_0))\in M$ and $\tm_\e,\mu_\e$ and $d_\e$ are defined in \eqref{mu} and \eqref{di}.
 According   this notation, we set  ${\boldsymbol { \omega}}={\boldsymbol { \omega}}(z)$ the function corresponding to ${\boldsymbol {\tilde\omega}}={\boldsymbol {\tilde\omega}}(y_0, y).$\\
Let $\delta>0$ be a fixed number with $4\delta<\hat{\delta}$, where $\hat{\delta}$ is given in \eqref{deta}. We consider a smooth cut-off function $\zeta_\delta(s)$ such that $\zeta_\delta(s)=1$ if $0<s<\delta$ and $\zeta_\delta(s)=0$ if $s>2\delta$. Let us consider the cut-off function $\eta^\e_\delta$ defined on the manifold $M_\rho$ by
$$\eta_\delta^\e(z)=\zeta_\delta\({{dist_g}(\xi,\Gamma)\over\rho}\) \qquad \mbox{for}\,\, \rho z=\xi\in M .$$   We remark that with this definition $\eta_\delta^\e(z)$ does not depend on the parameter functions.\\
We define our global first approximation of the problem \eqref{pb} ${\bf w}(z)$ as
\beq\label{appglobale}
{\bf w}(z)=\eta_\delta^\e(z){\boldsymbol { \omega}}(z).
\eeq
We look for a solution to problem \eqref{pbgluing} of the form $u={\bf w}+\Phi$,
namely
\beq\label{pbgluing2}
\Delta_g\Phi+p{\bf w}^{p-1}\Phi+N(\Phi)+E=0 \ \mbox{in}\   \M_\rho
\eeq
where
\begin{equation}\label{enne}N(\Phi)=\rho^{-\frac{N-2}{2}\e}({\bf w}+\Phi)^{p-\e}-{\bf w}^{p-\e}-p{\bf w}^{p-1}\Phi-\rho^2 h ({\bf w}+\Phi)\end{equation} and
\begin{equation}\label{errore-e} E=\Delta_g {\bf w}+{\bf w}^{p-\e}.\end{equation}
 We look for a solution $\Phi$ of \eqref{pbgluing2}   as $\Phi=\eta_{2\delta}\phi+\psi $
where the function $\phi$ is such that the corresponding function  $\tilde{\phi}$ via the change of variables \eqref{changeglu} is   defined only in $\D$.  It is immediate to check that $\Phi$ of this form solves \eqref{pbgluing2}  if the pair $(\psi, \phi)$ solves the following nonlinear coupled system:
\beq\label{gluing3}
\Delta_g \psi +(1-\eta_{2\delta}^\e)p{\bf w}^{p-1}\psi=-2\nabla_g \phi\nabla_g \eta_{2\delta}^\e-\phi\Delta_g \eta_{2\delta}^\e-(1-\eta_{2\delta}^\e)N(\eta_{2\delta}^\e \phi+\psi)\  \mbox{in}\  \M_\rho
\eeq
and
\beq\label{gluing2}
\A(\tilde{\phi})+p {\boldsymbol {\tilde\omega}}^{p-1}\tilde{\phi}= -{\mathcal N}(\zeta_{2\delta}^\e \tilde{\phi}+\tilde{\psi})-\S_\e({\boldsymbol {\tilde\omega}})-p{\boldsymbol {\tilde\omega}}^{p-1}\tilde\psi \  \mbox{in}\  \D,
\eeq
where
\begin{equation}\label{enne-fra}{\mathcal N} (\tilde \Phi)=\tilde\mu_\e^{-\frac{N-2}{2}\e}( {\boldsymbol {\tilde\omega}}+\tilde\Phi)^{p-\e}-{\bf w}^{p-\e}-p{ {\boldsymbol {\tilde\omega}}}^{p-1}\tilde\Phi-\tilde\mu_\e^2 \tilde h  \tilde\Phi ,\ \tilde \Phi=\zeta_{2\delta}^\e \tilde{\phi}+\tilde{\psi}.\end{equation}
Indeed,   problem \eqref{pbgluing2} in a scaled neighborhood of the geodesic looks like problem \ref{gluing2} and the error $E$ given in \eqref{errore-e} via the change of variables \eqref{changeglu} is nothing but the error term $ \S_\e({\boldsymbol {\tilde\omega}})$ defined in \eqref{experr2}.\\

Given $\phi$  such that $\tilde{\phi}$ is defined in $\D$, we first solve problem \eqref{gluing3} for $\psi$ (see Section 6 of  \cite{DMP}).
\begin{lemma}\label{lem-glu1}
For any $R>0$ there exists $r>0$ such that for any function $\phi$   such that the corresponding function  $\tilde{\phi}$   is   defined only in $\D$ with  $\|\tilde{\phi}\|_*\le r$, there exists a unique solution  $\psi=\psi(\phi)$ of \eqref{gluing3} with
$$\|\psi\|_\infty\leq R \e^{\frac{N-4}{2}}\|\tilde{\phi}\|_*.$$ Moreover, the nonlinear operator $\psi$ satisfies a Lipschitz condition of the form
\beq\label{gluing4}
\|\psi(\phi_1)-\psi(\phi_2)\|_\infty\leq c \e^{\frac{N-4}{2}}\|\phi_1-\phi_2\|_*,
\eeq
for some positive constant $c$ independent on $\e.$
\end{lemma}

Finally, we substitute $\tilde{\psi}=\tilde{\psi}( \phi )$  (via the change of variables \eqref{changeglu}) in the equation \eqref{gluing3} and we reduce the full problem \eqref{pb} to solving the following (nonlocal) problem in $\D$:
\beq\label{pbgluingdef}
\A(\tilde{\phi})+p{\boldsymbol {\tilde\omega}}^{p-1}\tilde{\phi}=-{\mathcal N}(\eta_{2\delta}^\e\tilde{\phi}+\tilde{\psi}( {\phi}))-\S_\e({\boldsymbol {\tilde\omega}})-p{\boldsymbol {\tilde\omega}}^{p-1}\tilde{\psi}( {\phi})  \  \mbox{in}\  \D.
\eeq

\subsection{The nonlinear projected problem}\label{nonlineare}
We can  solve the following projected problem associated to \eqref{pbgluingdef}: {\it given $\mu, d$ and $e$ satisfying \eqref{norma},   find functions $\tilde \phi$ and $c_j(y_0)$ for $j=0, \ldots, N+1$ such that}
\beq\label{non2}
\left\{
\begin{array}{lr}
L(\tilde \phi)=-S_\e({\boldsymbol {\tilde\omega}})+ {\mathfrak N}(\tilde\phi)+ \sum_{j=0}^N c_j Z_j\qquad \mbox{in}\,\,\, \D \\\\
\tilde\phi\(y_0+{2\ell\over\rho}, y\)= \phi(y_0, Ay) \qquad \mbox{for any}\,\,\, (y_0, y)\in \D,\\\\
\int_{\D_{y_0}}\tilde\phi Z_j  dy =0 \ \mbox{and for any}\  y_0\in\[-{\ell\over\rho},{\ell\over\rho}\], \ j=0,1, \ldots, N+1.
\end{array}
\right.
\eeq
   Here $S_\e({\boldsymbol {\tilde\omega}})$ is given in \eqref{experr2} and
    \begin{align*}
& L(\tilde\phi):= \A (\tilde\phi)+p \omega^{p-1}\tilde\phi \ \hbox{($ \A $ is    in Lemma \ref{lemmalaplace} and   $\omega$ is  in \eqref{ansa})},\\
& {\mathfrak N}(\tilde\phi):=p(\omega^{p-1}-{\boldsymbol{\tilde \omega}}^{p-1})\tilde\phi- {\mathcal N}(\zeta^{\e}_{2\delta}\tilde\phi+\tilde\psi( \phi))-p{\boldsymbol{\tilde \omega}}^{p-1}\tilde\psi( \phi)\ \hbox{($ {\mathcal N} $   is    in \eqref{enne-fra})}.\end{align*}

\begin{proposition}\label{propfixed}
There exists $c>0$ such that for all sufficiently small $\e$ and all $\mu, d$ and $e$ satisfying \eqref{norma}, problem \eqref{non2} has a unique solution $\tilde\phi=\tilde\phi(\mu, d, e)$ and $c_j=c_j(\mu, d, e)$ which satisfies
\beq\label{stimaphi}
\|\phi\|_*\leq c\e^{\frac 32}.
\eeq
Moreover, $\tilde\phi$ depends Lipschitz continuously on $\mu,$ $d$ and $e$ in the sense
$$\|\tilde\phi(\mu_1,d_1,e_1)-\tilde\phi(\mu_2,d_2,e_2)\|_*\le \e^{5\over2}\|(\mu_1-\mu_2,d_1-d_2,e_1-e_2)\|$$
for some  positive constant $c$ independent of $\e$ and uniformly with respect to $\mu, d$ and $e$ which satisfy \eqref{norma}.

\end{proposition}
\begin{proof} We argue exactly as in Section 7 of \cite{DMP}, using a contraction mapping argument and the linear theory developed in Proposition \ref{inv}.
 \end{proof}

\section{The reduced problem}\label{equazioni}
\subsection{The reduced  system}
We find $N+1$ equations relating $\mu, d$ and $e$ to get all the coefficients $c_j$ in \eqref{non2} identically equal to zero. To do this, we multiply equation \eqref{non2} by $Z_j$, for all $j=0, \ldots, N+1$ and we integrate in $y$. Thus, the system $$c_j(\rho y_0)=0, \ j=0,1, \ldots, N+1$$ is equivalent to
$$\int_{\D_{y_0}}S_\e({\boldsymbol {\tilde\omega}})Z_j\, dy +\int_{\D_{y_0}}\left(L(\tilde\phi)-{\mathfrak N}(\tilde\phi)\right)Z_j\, dy =0, \ j=0,1, \ldots, N+1,$$ for any  $y_0\in\[-{\ell\over\rho},{\ell\over\rho}\] .$\\
By Proposition \ref{propfixed} it follows that
$$\int_{\D_{y_0}}\left(L(\tilde\phi)-{\mathfrak N}(\tilde\phi)\right)Z_j\, dy =\e^3 \theta,$$ where $\theta=\theta(\rho y_0)$ is as in Lemma \ref{errore2}.

Hence the equations $c_j=0$ are equivalent to the following limit system on $N+2$ nonlinear ordinary differential equations:
\beq\label{sistemaequazioni}
\left\{
\begin{array}{lr}
L_{N+1}(\mu):=- \ddot{\mu}+\(a_n\sigma \pm{b_n\over\mu^2_0} \) \mu =-\alpha_{N+1}(x_0)-c_3Q(x_0,d)+\e |\ln \e|M_{N+1}\\\\
L_k(d ):=-\ddot{d}_k+\sum_{j=1}^N R_{0j0k}d_j= \sqrt\e M_k,\  k=1, \ldots, N\\\\
L_0(e):=\e a_0 \ddot{e} +\lambda_1 e =-\alpha_0(x_0)-  c_4Q(x_0,d) -\beta(x_0) \mu+\e  |\ln\e|M_0
\end{array}
\right.
\eeq
where $\mu, d_1, \ldots, d_N,e\in C^2_{2\ell}(\mathbb R)$ and
\begin{itemize}
\item[-] the functions $\alpha_i$ and $\beta$ are explicit functions of $x_0$, smooth and uniformly bounded in $\e$ given in Lemma \ref{errore2}
\item[-] the operator $Q $ is quadratic in $d$ (see Lemma \ref{errore2}) and it is uniformly bounded in $L^\infty_{2\ell}(\mathbb R)$ for $(\mu, d, e)$ satisfying \eqref{norma}
\item[-] the operators $M_i=M_i(\mu,d,e)$ can be decomposed as $M_i(\mu,d,e)=A_i(\mu,d,e)+K_i(\mu,d,e)$ where
\begin{itemize}\item[-] $K_i$ is uniformly bounded in $L^\infty_{2\ell}(\mathbb R)$ for $(\mu, d, e)$ satisfying \eqref{norma} and it is compact
\item[-] $A_i$ depends on $(\mu,d,e)$ and their first and second derivatives and it satisfies
$$\|A_i(\mu_2,d_2,e_2)-A_i(\mu_1,d_1,e_1)\|\le o(1)\|(\mu_2-\mu_1,d_2-d_1,e_2-e_1)\|$$
uniformly  for $(\mu, d, e)$ satisfying \eqref{norma}
\item[-] the dependance on  $(\ddot\mu,\ddot d,\ddot e)$ is linear
\end{itemize}
\end{itemize}

Our goal is to solve \eqref{sistemaequazioni} in $\mu, d$ and $e$.
To do so, we first analyze the invertibility of the linear operator  $L_{N+1}$.\\
\begin{lemma}\label{lnp}
For any $f\in L^\infty_{2\ell}(\mathbb R),$   there exists a unique $\mu\in C^2_{2\ell}(\mathbb R)$ solution of
 $L_{N+1}(\mu)=f .$
Moreover, there exists $c$ such that
$$\|\mu\|_\infty+\|\dot{\mu}\|_\infty\leq c \|f\|_\infty.$$
\end{lemma}
\begin{proof}
The non degeneracy condition of the solution $\mu_0$  translates into the fact that the periodic O.D.E.
$$- \ddot{\mu}+\(a_n\sigma \pm{b_n\over\mu^2_0} \) \mu=0\ \hbox{in}\ [0,2\ell]$$
has only the trivial solutions. Therefore the claim follows.
\end{proof}
Next, we  analyze the invertibility of the linear operator  $L_{0}$.\\
\begin{lemma}\label{l0}
Assume
$$|\e m^2-\kappa^2|>\nu \sqrt{\e}\ \hbox{for any}\ m=1, 2, \ldots
$$
for some $\nu$ positive, where
$$\kappa:={\pi \over2}\sqrt{\lambda_1}\int_{-\ell}^{+\ell}{1\over\sqrt{a_0(s)}}ds.$$
 For any $f\in C^0_{2\ell}(\mathbb R)\cap L^\infty_{2\ell}(\mathbb R),$  there exists a unique solution $e\in C^2_{2\ell}(\mathbb R)$ of $L_0(e)=f.$   Moreover, there exists $c$ such that
$$\e\|\ddot{e}\|_\infty+\sqrt\e\|\dot{e}\|_\infty +\|e\|_\infty\leq c{1\over \sqrt\e}\|f\|_\infty,$$
Finally, if $f\in C^2_{2\ell}(\mathbb R)$, then
$$\e\|\ddot{e}\|_\infty+\sqrt\e\|\dot{e}\|_\infty+\|e\|_\infty\leq c\left[\|\ddot{f}\|_\infty+\|\dot{f}\|_\infty+\|f\|_\infty\right].$$
\end{lemma}
\begin{proof}
We argue as in in Lemma 8.2 of \cite{DMP}.
\end{proof}
Finally, we consider the invertibility of the linear operator $(L_1,\dots,L_N).$
\begin{lemma}\label{lk}
Assume the geodesic is non degenerate.
For any $f=(f_1,\dots,f_N)$ with $f_k\in L^\infty_{2\ell}(\mathbb R)$, there exists a $ d=(d_1,\dots,d_N)$ with $d_k\in C^2_{2\ell}(\mathbb R)$  such that $L_k(d)=f_k$ for any $k=1,\dots,N.$
 Moreover, there exists $c$ such that
 $$\|\ddot{d} \|_\infty+\|\dot{d} \|_\infty+\|d \|_\infty \leq c\|f\|_\infty.$$
\end{lemma}
\begin{proof}
It is useful to point out  that assumption \eqref{jacobi} about non degeneracy of $\Gamma$ in normal coordinates translates exactly into the fact that the linear system of O.D.E.'s
$$-\ddot{ {d}}_k+\sum_{j=1}^N R_{0j0k} {d}_j=0,\ \hbox{in}\ [0,2\ell],\quad     k=1, \ldots, N,\
$$
has only the trivial solution $ {d}\equiv 0$  satisfying the periodicity condition \eqref{periodicA}. Therefore, the claim  follows.\\
\end{proof}
\subsection{The choice of parameters: the proof completed!}
Now, we are ready to complete the proof, finding parameters which solve the reduced problem \eqref{sistemaequazioni}.\\\\
First,   by Lemma \ref{lnp} we find
  $\hat\mu_0$    solution  of  $$L_{N+1}(\hat\mu_0)=-\alpha_{N+1}(x_0) ,\ \hbox{ with}\
   \|\ddot{\hat\mu}_0\||_\infty+\|\dot{\hat\mu}_0\|_\infty+\|\hat\mu_0\|_\infty\leq c.$$ Then, by Lemma \ref{l0} we find $\hat e_0$ solution of
$$L_0(\hat{e}_0)=-\alpha_0 -\beta\hat\mu_0, ,\ \hbox{ with}\ \e \|\ddot{\hat e}_0\|_\infty+\sqrt\e\|\dot{\hat e}_0\|_\infty+\|\hat e_0\|_\infty\leq c.$$
Therefore,  $\|(\hat\mu_0,0,\hat e_0)\|\leq c.$  Let us define $$ \mu=\hat\mu_0+ \hat\mu_1,\ d=   \hat{d}_1,\  e=\hat e_0+ \hat e_1.$$ The system \eqref{sistemaequazioni} reduces to
\beq\label{sistemaequazioniridotto}
\left\{
\begin{array}{lr}
L_{N+1}(\hat\mu_1)=-c_3Q(x_0,\hat{d}_1)+\e|\ln\e|   M_{N+1} \\\\
L_k(\hat{d}_1)=\sqrt\e M_k ,\ k=1,\dots,N\\\\
L_0(\hat e_1)=-c_4Q(x_0,\hat{d}_1)-\beta(x_0)\hat\mu_1+\e|\ln\e|    M_0
\end{array}
\right.
\eeq
Let us observe now that the linear operator $$\mathcal L(\hat\mu_1, \hat d_1, \hat e_1)= (L_{N+1}(\hat \mu_1), L_N(\hat d_{1  }), \ldots, L_1(\hat d_{1 }), L_0(\hat e_1))$$ is invertible with bounds for $\mathcal L(\hat \mu_1,\hat  d_1,\hat  e_1)=(f, g, h)$ given by $$\|(\hat \mu_1,\hat  d_1,\hat  e_1)\|\leq C\left[\|f\|_\infty+\|g\|_\infty+\e^{-1/2}\|h\|_\infty\right].$$ Finally, by the contraction mapping principle it follows that,   the problem
  \eqref{sistemaequazioniridotto}
has a unique solution  with $$\|\hat\mu_1\|_\infty <c\e|\ln\e|,\qquad \|\hat d_1\|_\infty<\sqrt\e, \qquad \|\hat e_1\|_\infty<\sqrt\e|\ln\e|.$$ That concludes the proof.

\section{The linear theory}\label{lineare}
Here we recall  a linear theory necessary to solve problem \eqref{pA}, which has been developed in Section 3 of \cite{DMP}.\\
 Let us consider the operator  $\mathcal L_0:=\Delta_{\mathbb R^N}+pw^{p-1} $. It is well-known that the $L^2-$ null space of the operator $\mathcal L_0$ is $N+1-$ dimensional and spanned by the functions $$Z_j(y):=\partial_{j}w(y),\  j=1, \ldots, N\ \hbox{and}\ Z_{N+1}(y):=y\cdot\nabla w(y)+\frac{N-2}{2}w(y).$$ Moreover it is known that (see \cite{DMP}) that the operator $\mathcal L_0$ has one negative eigenvalue $-\lambda_1<0,$ whose corresponding eigenfunction $Z_0$    (normalized to have $L^2-$ norm equal to $1$) decays exponentially at infinity with exponential order $O(e^{-\sqrt{\lambda_1}|x|})$.\\\\ The following results (see Lemma 3.1 of \cite{DMP} and also \cite{DPV}) are useful in order to obtain a priori estimates and a solvability theory for problem \eqref{pA}.
\begin{lemma}\label{lemmalin1}
Assume that $\lambda\not\in\left\{0, \pm\sqrt{\lambda_1}\right\}$. Then for $g\in L^\infty(\mathbb R^N)$, there exists a unique bounded solution of $$(\mathcal L_0-|\lambda|^2)\psi=g$$ in $\mathbb R^N$. Moreover $$\|\psi\|_{L^\infty}\leq c_{\lambda}\|g\|_{L^\infty}$$ for some constant $c_{\lambda}>0$ only depending on $\lambda$.
\end{lemma}
 \begin{lemma}\label{lemmalin2}
Let $\phi$ a bounded solution of $$\partial_{00}\phi+\Delta_y\phi+p w^{p-1}\phi=0\qquad \mbox{in}\qquad \mathbb R^{N+1}.$$ Then $\phi(y_0, y)$ is a linear combination of the functions $Z_j$, $j=1, \ldots, N+1$, $Z_0(y)\cos(\sqrt{\lambda_1}y_0)$, $Z_0(y)\sin(\sqrt{\lambda_1}y_0)$.
\end{lemma}
 Now, we study   a slightly more general problem than \eqref{pA} that involves the essential features needed. For any constant $M>0$ we consider the domain $\D$ defined as
\beq\label{D}
\D:=\left\{(y_0, y)\in\mathbb R\times \mathbb R^{N }\quad:\quad  |y|<M\right\}
\eeq
and  given     a function $\phi$ defined on $\D$, an operator of the form
$$L(\phi):=b(y_0)\partial_{00}\phi+\Delta_y \phi+pw^{p-1}\phi+\sum\limits_{i,j }  b_{ij}(y_0, y)\partial_{ij}\phi+\sum\limits_{i  }  b_i(y_0, y)\partial_i\phi+d(y_0, y)\phi .$$    Then for a given function $g$ we want to solve the following projected problem:
\beq\label{lin1}
\left\{
\begin{array}{lr}
L(\phi)=g +\sum_{j=0}^{N+1}c_j(y_0)Z_j(y)\qquad \mbox{in}\,\ \D\\\\
\int_{\D_{y_0}}\phi(y_0, y)Z_j(y)\, dy=0\qquad \mbox{for any}\  y_0\in \mathbb R, \  j=0, \ldots, N
\end{array}
\right.
\eeq
where $$\D_{y_0}:=\left\{y\in \mathbb R^N : (y_0, y)\in \D\right\}.$$
We fix a number $2\leq \nu <N$ and consider the $L^\infty-$ weighted norms
\begin{align*}
\|\phi\|_* & :=\sup_{\D}(1+|y|^{\nu-2})|\phi(y_0, y)|+\sup_{\D}(1+|x|^{\nu-1})|D\phi(x_0, x)| ,
\\
\|g\|_{**} & :=\sup_{\D}(1+|y|^\nu)|g(y_0, y)|.
\end{align*}
We assume that all functions involved are smooth. The following result (see Proposition 3.2 of \cite{DMP}) establishes  existence and uniform a priori estimates for problem \eqref{lin1} in the above norms, provided that appropriate bounds for the coefficients hold.\\
\begin{proposition}\label{inv}
Assume that $N\geq 7$ and $N-2\leq \nu <N$. Assume that there exists   $m>0$ such that $$m\leq b(y_0)\leq m^{-1}\qquad \mbox{for any}\  y_0\in \mathbb R.$$ There exist $\delta>0$ and $C>0$ such that if
\beq\label{inv1}
 M\|\partial_0 b\|_\infty+\sum\limits_{i,j }  \(\|b_{ij}\|_\infty +\|D b_{ij}\|_\infty\)+\sum\limits_{i  }  \|(1+|y|)b_i\|_\infty+\|(1+|y|^2)d\|_\infty<\delta
\eeq
then for any $g$ with $\|g\|_{**}<\infty$ there exists a unique solution $\phi=T(g)$ of problem \eqref{lin1} with $\|\phi\|_*<\infty$ and  it holds true that $$\|\phi\|_*\leq C\|g\|_{**}.$$
\end{proposition}

\section{Appendix}

\subsection{Proof of \eqref{metrica}}
Let $E_0,E_1,\dots,E_N$ the coordinate vectors as given in the Introduction.
 By our choice of coordinates it follows that $\nabla_E E=0$ on $\Gamma$ for any vector field $E$, that is a linear combination (with coefficients depending only on $x_0$) of the $E_j$'s, $j=1, \ldots, N$.\\ In particular, for any $i, j=1, \ldots, N$ and for any $t\in\mathbb R$, we have $\nabla _{E_i+t E_j}(E_i+t E_j)=0$ on $\Gamma$, which implies $\nabla_{E_i}E_j+\nabla_{E_j}E_i=0$ for every $i, j=1, \ldots, N$.\\ Using the fact that $E_i$'s are coordinate vectors for $j=1, \ldots, N$ and in particular $\nabla_{E_a}E_b=\nabla_{E_b}E_a$ for all $a, b=0, \ldots, N$, we obtain that $\nabla E_j E_i=0$ for every $i, j=1, \ldots, N$. The geodesic coordinate for $\Gamma$ translates precisely into $\nabla E_0 E_0=0$.\\ These facts immediately yields
\begin{equation}\label{+secondo}
\partial_m g_{ij}=E_m\langle E_i, E_j\rangle=\langle \nabla_{E_m}E_i, E_j\rangle+\langle E_i, \nabla_{E_m}E_j\rangle=0
\end{equation}
on $\Gamma$ with $i, j, m= 1, \ldots, N$.\\\\Moreover, since $E_a$'s are coordinate vectors for $a=0, \ldots, N$, we obtain
\begin{eqnarray}\label{+}
\nonumber
\partial_m g_{0j}&=& E_m \langle E_0, E_j\rangle\\
\nonumber
&=& \langle\nabla_{E_m}E_0, E_j\rangle+\langle E_0, \nabla_{E_m}E_j\rangle\\
&=& \langle \nabla_{E_0} E_m, E_j\rangle+\langle E_0, \nabla_{E_m}E_j\rangle =0
\end{eqnarray}
on $\Gamma$ with $m, j=1, \ldots, N$.\\\\ Here we used the fact that $\nabla_{E_0}E_m=0$ on $\Gamma$, namely that $\nabla_{E_0}E_m$ has zero normal components.\\ Moreover by \eqref{+secondo} it follows that
\begin{equation}\label{+primo}
\partial_m g_{00}=0\qquad \mbox{on}\quad \Gamma.
\end{equation}
We can also prove that the components $R_{0m0j}$ of the curvature tensor are given by
\begin{equation}\label{++}
R_{0m0j}=-\frac 12 \partial_{mj} g_{00}.
\end{equation}
Indeed, we have
\begin{eqnarray*}
-R_{0m0j}&=&\langle R(E_0, E_j)E_0, E_m\rangle\\
&=& \langle \nabla_{E_0}E_j E_0, E_m\rangle-\langle \nabla_{E_j}\nabla_{E_0} E_0, E_m\rangle\\
&=& \langle \nabla_{E_0}\nabla E_j E_0, E_m\rangle -E_j \langle \nabla_{E_0} E_0, E_m\rangle-\langle \nabla_{E_0}E_0, \nabla_{E_j}E_m\rangle\\
&=& \langle \nabla_{E_0}\nabla_{E_j}E_0, E_m\rangle-E_j\langle \nabla_{E_0}E_0, E_m\rangle\\
&=& \langle \nabla_{E_0}\nabla_{E_j}E_0, E_m\rangle-E_j E_0\rangle E_0, E_m\rangle+E_j\langle E_0, \nabla_{E_0}E_m\rangle\\
&=& \langle \nabla_{E_0}\nabla_{E_j}E_0, E_m\rangle + E_j\langle E_0, \nabla_{E_m}E_0\rangle\\
&=& \frac 12 E_j E_m \langle E_0, E_0\rangle + E_0\langle \nabla_{E_j}E_0, E_m\rangle-\langle\nabla_{E_j}E_0,\nabla_{E_0}E_m\rangle\\
&=&\frac 12 \partial_{mj} g_{00}
\end{eqnarray*}
where here we have used the above properties and the fact that
$$\nabla_{E_j}E_0=\nabla_{E_0}E_j=\frac 12 \partial_j g_{00}E_0=0.$$ By \eqref{+}, \eqref{++}, \eqref{+primo} and \eqref{+secondo} the claim follows.

\end{document}